# PDES FOR THE JOINT DISTRIBUTIONS OF THE DYSON, AIRY AND SINE PROCESSES

By Mark Adler[1] and Pierre van Moerbeke[2]

*Brandeis University and Université de Louvain*

In a celebrated paper, Dyson shows that the spectrum of an $n \times n$ random Hermitian matrix, diffusing according to an Ornstein–Uhlenbeck process, evolves as $n$ noncolliding Brownian motions held together by a drift term. The universal edge and bulk scalings for Hermitian random matrices, applied to the Dyson process, lead to the Airy and Sine processes. In particular, the Airy process is a continuous stationary process, describing the motion of the outermost particle of the Dyson Brownian motion, when the number of particles gets large, with space and time appropriately rescaled.

In this paper, we answer a question posed by Kurt Johansson, to find a PDE for the joint distribution of the Airy process at two different times. Similarly we find a PDE satisfied by the joint distribution of the Sine process. This hinges on finding a PDE for the joint distribution of the Dyson process, which itself is based on the joint probability of the eigenvalues for coupled Gaussian Hermitian matrices. The PDE for the Dyson process is then subjected to an asymptotic analysis, consistent with the edge and bulk rescalings. The PDEs enable one to compute the asymptotic behavior of the joint distribution and the correlation for these processes at different times $t_1$ and $t_2$, when $t_2 - t_1 \to \infty$, as illustrated in this paper for the Airy process. This paper also contains a rigorous proof that the extended Hermite kernel, governing the joint probabilities for the Dyson process, converges to the extended Airy and Sine kernels after the appropriate rescalings.

Received August 2003; revised October 2004.

[1]Supported by NSF Grant DMS-01-00782.

[2]Supported by NSF Grant DMS-01-00782, and grants from European Science Foundation, NATO, FNRS and Francqui Foundation. This work was done while the author was a member of the Clay Mathematics Institute, One Bow Street, Cambridge, Massachusetts 02138, USA.

**1. Stating the results.** The *Dyson Brownian motion* [4]

$$(\lambda_1(t), \ldots, \lambda_n(t)) \in \mathbb{R}^n$$

with transition density $p(t, \mu, \lambda)$ satisfies the diffusion equation

$$\frac{\partial p}{\partial t} = \frac{1}{2} \sum_1^n \frac{\partial}{\partial \lambda_i} \Phi(\lambda) \frac{\partial}{\partial \lambda_i} \frac{1}{\Phi(\lambda)} p$$

$$= \sum_1^n \left( \frac{1}{2} \frac{\partial^2}{\partial \lambda_i^2} - \frac{\partial}{\partial \lambda_i} \frac{\partial \log \sqrt{\Phi(\lambda)}}{\partial \lambda_i} \right) p$$

with

$$\Phi(\lambda) = \Delta^2(\lambda) \prod_1^n e^{-\lambda_i^2/a^2}.$$

Roughly speaking, it represents $n$ Brownian motions repelling one another, with the exponential in $\Phi(\lambda)$ having the effect of preventing the system from flying out to infinity. In his beautiful paper, Dyson generalizes the random matrix ensembles in such a way that *the Coulomb gas model acquires a meaning as a dynamical system, rather than a static model*. He shows the repelling Brownian motion above corresponds to the motion of the eigenvalues $(\lambda_1(t), \ldots, \lambda_n(t))$ of a Hermitian matrix $B$, evolving according to the Ornstein–Uhlenbeck process

$$(1.1) \qquad \frac{\partial P}{\partial t} = \sum_{i,j=1}^{n^2} \left( \frac{1}{4}(1+\delta_{ij}) \frac{\partial^2}{\partial B_{ij}^2} + \frac{1}{a^2} \frac{\partial}{\partial B_{ij}} B_{ij} \right) P,$$

with transition density $(c := e^{-t/a^2})$

$$P(t, \bar{B}, B) = Z^{-1} \frac{1}{(1-c^2)^{n^2/2}} e^{-(1/(a^2(1-c^2))) \operatorname{Tr}(B - c\bar{B})^2}.$$

The $B_{ij}$'s in (1.1) denote the $n^2$ free (real) parameters in the Hermitian matrix $B$, with the $B_{ii}$'s being its diagonal elements. In the limit $t \to \infty$, this distribution tends to the stationary distribution

$$Z^{-1} e^{-(1/a^2) \operatorname{Tr} B^2} dB = Z^{-1} \Delta^2(\lambda) \prod_1^n e^{-\lambda_i^2/a^2} d\lambda_i.$$

With this invariant measure as initial condition, the joint distribution reads

$$(1.2) \quad \begin{aligned} &P(B(0) \in dB_1, B(t) \in dB_2) \\ &= Z^{-1} \frac{dB_1 dB_2}{(1-c^2)^{n^2/2}} e^{-(1/(a^2(1-c^2))) \operatorname{Tr}(B_1^2 - 2cB_1 B_2 + B_2^2)}, \end{aligned}$$

for which a nonlinear PDE will be found in Theorem 1.1. According to [9], the joint probabilities for the $\lambda_i$'s can also be expressed in terms of a Fredholm determinant of the so-called extended Hermite kernel $\hat{K}^{H,n}_{t_i t_j}(x,y)$, a matrix kernel, defined in Section 7. As elaborated in that section, for $E_1$ and $E_2 \subset \mathbb{R}$, we have, for $1 \leq i \leq n$,

$$(1.3)\ P(\text{all } \lambda_i(t_1) \in E_1,\ \text{all } \lambda_i(t_2) \in E_2) = \det(I - (\chi_{E_k^c} K^{H,n}_{t_k t_\ell} \chi_{E_\ell^c})_{1 \leq k,\ell \leq 2}).$$

Since expression (1.2) is symmetric in $B_1$ and $B_2$, the probability (1.3) for the Dyson process is symmetric in $E_1$ and $E_2$. Throughout the paper, we normalize the problem, by setting $a = 1$.

The *Airy process* is defined by an appropriate rescaling of the largest eigenvalue $\lambda_n$ in the Dyson diffusion,

$$(1.4) \qquad A(t) = \lim_{n \to \infty} \sqrt{2} n^{1/6} (\lambda_n(n^{-1/3}t) - \sqrt{2n}),$$

in the sense of convergence of distributions for a finite number of $t$'s. This process was introduced by Prähofer and Spohn [12] in the context of polynuclear growth models and further investigated by Johansson [8]. Prähofer and Spohn showed the Airy process is a stationary process with continuous sample paths; thus the probability $P(A(t) \leq u)$ is independent of $t$, and is given by the Tracy–Widom distribution [14],

$$(1.5) \qquad P(A(t) \leq u) = F_2(u) := \exp\left(-\int_u^\infty (\alpha - u) q^2(\alpha)\, d\alpha\right),$$

with $q(\alpha)$ the solution of the *Painlevé* II equation,

$$(1.6) \quad q'' = \alpha q + 2q^3 \quad \text{with } q(\alpha) \cong \begin{cases} -\dfrac{e^{-(2/3)\alpha^{3/2}}}{2\sqrt{\pi}\alpha^{1/4}}, & \text{for } \alpha \nearrow \infty, \\ \sqrt{-\alpha/2}, & \text{for } \alpha \searrow -\infty. \end{cases}$$

Here, the joint probabilities for the process $A(t)$ can also be expressed in terms of the Fredholm determinant of the extended Airy kernel $\hat{K}^A_{t_i t_j}(x,y)$, another matrix kernel, which is an appropriate limit of the extended Hermite kernel above; see [6, 8, 11, 12]. It leads to

$$(1.7) \qquad P(A(t_1) \in E_1, A(t_2) \in E_2) = \det(I - (\chi_{E_i^c} K^A_{t_i t_j} \chi_{E_j^c})_{1 \leq i,j \leq 2}),$$

which is also symmetric in $E_1$ and $E_2$, as a consequence of the symmetry for the Dyson process.

At MSRI (Sept. 2002), Kurt Johansson, whom we thank for introducing us to the Airy process, posed the question whether a PDE can be found for the joint probability of this process; see [8]. The present paper answers this question (Theorem 1.2), which enables us to derive the asymptotics of the large time correlations for the Airy process (Theorem 1.6); this question was posed in [12]. Our results on the Airy process for the special case of semi-infinite intervals appeared in [2], as well as the asymptotics.

The *Sine process* is an infinite collection of noncolliding processes $S_i(t)$, obtained by rescaling the bulk of the Dyson process in the same way as the bulk of the spectrum of a large Gaussian random matrix; namely,

$$(1.8) \qquad S_i(t) := \lim_{n\to\infty} \frac{\sqrt{2n}}{\pi} \lambda_{n/2+i}\left(\frac{\pi^2 t}{2n}\right) \qquad \text{for } -\infty < i < \infty$$

in the sense of convergence of distributions for a finite number of $t$'s. This process was defined by Tracy and Widom in [17]. Similarly, by taking the bulk scaling limit of the extended Hermite kernel (1.3), the joint probabilities for the $S_i$'s can also be expressed in terms of the Fredholm determinant of the extended sine kernel $K^S_{t_i t_j}(x,y)$, yet another matrix kernel,

$$(1.9)\, P(\text{all } S_i(t_1) \in E_1^c,\ \text{all } S_i(t_2) \in E_2^c) = \det(I - (\chi_{E_k} K^S_{t_k t_\ell} \chi_{E_\ell})_{1\leq k,\ell \leq 2}),$$

where here $E_1$ and $E_2$ must be compact for it to make sense. Note this probability is, as usual, symmetric in $E_1$ and $E_2$. These kernels will be discussed in Section 7. For this process, Theorem 1.4 gives a PDE for the joint probabilities.

The disjoint union of intervals

$$E_1 := \bigcup_{i=1}^r [a_{2i-1}, a_{2i}] \quad \text{and} \quad E_2 := \bigcup_{i=1}^s [b_{2i-1}, b_{2i}] \subseteq \mathbb{R},$$

and $t = t_2 - t_1$ specify linear operators, setting $c = e^{-t}$,

$$(1.10) \quad \begin{aligned}
\mathcal{A}_1 &= \sum_1^{2r} \frac{\partial}{\partial a_j} + c \sum_1^{2s} \frac{\partial}{\partial b_j}, \\
\mathcal{B}_1 &= c \sum_1^{2r} \frac{\partial}{\partial a_j} + \sum_1^{2s} \frac{\partial}{\partial b_j}, \\
\mathcal{A}_2 &= \sum_1^{2r} a_j \frac{\partial}{\partial a_j} + c^2 \sum_1^{2s} b_j \frac{\partial}{\partial b_j} + (1-c^2)\frac{\partial}{\partial t} - c^2, \\
\mathcal{B}_2 &= c^2 \sum_1^{2r} a_j \frac{\partial}{\partial a_j} + \sum_1^{2s} b_j \frac{\partial}{\partial b_j} + (1-c^2)\frac{\partial}{\partial t} - c^2.
\end{aligned}$$

The duality $a_i \leftrightarrow b_j$ reflects itself in the duality $\mathcal{A}_i \leftrightarrow \mathcal{B}_i$. We now state:

THEOREM 1.1 (Dyson process). *Given $t_1 < t_2$ and $t = t_2 - t_1$, the logarithm of the joint distribution for the Dyson Brownian motion $(\lambda_1(t), \ldots, \lambda_n(t))$,*

$$(1.11) \quad G_n(t; a_1, \ldots, a_{2r}; b_1, \ldots, b_{2s}) := \log P(\text{all } \lambda_i(t_1) \in E_1,\ \text{all } \lambda_i(t_2) \in E_2),$$

*satisfies a third-order nonlinear PDE in the boundary points of $E_1$ and $E_2$ and $t$, which takes on the simple form, setting $c = e^{-t}$,*

$$(1.12) \qquad \mathcal{A}_1 \frac{\mathcal{B}_2 \mathcal{A}_1 G_n}{\mathcal{B}_1 \mathcal{A}_1 G_n + 2nc} = \mathcal{B}_1 \frac{\mathcal{A}_2 \mathcal{B}_1 G_n}{\mathcal{A}_1 \mathcal{B}_1 G_n + 2nc}.$$

The proof of this theorem will be given in Section 3.

Similarly, the disjoint union of intervals

$$
(1.13) \quad E_1 := \bigcup_{i=1}^{r} [u_{2i-1}, u_{2i}] \quad \text{and} \quad E_2 := \bigcup_{i=1}^{s} [v_{2i-1}, v_{2i}] \subseteq \mathbb{R},
$$

and $t = t_2 - t_1$ define another set of linear operators

$$
(1.14) \quad \begin{aligned} L_u &:= \sum_{1}^{2r} \frac{\partial}{\partial u_i}, & L_v &:= \sum_{1}^{2s} \frac{\partial}{\partial v_i}, \\ E_u &:= \sum_{1}^{2r} u_i \frac{\partial}{\partial u_i} + t \frac{\partial}{\partial t}, & E_v &:= \sum_{1}^{2s} v_i \frac{\partial}{\partial v_i} + t \frac{\partial}{\partial t}. \end{aligned}
$$

We now give the equations for the joint probabilities of the Airy and Sine processes, which will be proven in Section 4:

THEOREM 1.2 (Airy process). *Given $t_1 < t_2$ and $t = t_2 - t_1$, the joint distribution for the Airy process $A(t)$,*

$$
G(t; u_1, \ldots, u_{2r}; v_1, \ldots, v_{2s}) := \log P(A(t_1) \in E_1, A(t_2) \in E_2),
$$

*satisfies a third-order nonlinear PDE in the $u_i$, $v_i$ and $t$, in terms of the Wronskian $\{f(y), g(y)\}_y := f'(y) g(y) - f(y) g'(y)$,*

$$
(1.15) \quad \begin{aligned} &((L_u + L_v)(L_u E_v - L_v E_u) + t^2 (L_u - L_v) L_u L_v) G \\ &= \tfrac{1}{2} \{(L_u^2 - L_v^2) G, (L_u + L_v)^2 G\}_{L_u + L_v}. \end{aligned}
$$

COROLLARY 1.3. *In the case of semi-infinite intervals $E_1$ and $E_2$, the PDE for the Airy joint probability,*

$$
H(t; x, y) := \log P\left( A(t_1) \leq \frac{y+x}{2}, A(t_2) \leq \frac{y-x}{2} \right),
$$

*takes on the following simple form in $x, y$ and $t^2$, with $t = t_2 - t_1$, also in terms of the Wronskian,*

$$
(1.16) \quad 2t \frac{\partial^3 H}{\partial t \, \partial x \, \partial y} = \left( t^2 \frac{\partial}{\partial x} - x \frac{\partial}{\partial y} \right) \left( \frac{\partial^2 H}{\partial x^2} - \frac{\partial^2 H}{\partial y^2} \right) + 8 \left\{ \frac{\partial^2 H}{\partial x \, \partial y}, \frac{\partial^2 H}{\partial y^2} \right\}_y.
$$

REMARK. Note for the solution $H(t; x, y)$,

$$
\lim_{t \searrow 0} H(t; x, y) = \log F_2 \left( \min\left( \frac{y+x}{2}, \frac{y-x}{2} \right) \right).
$$

THEOREM 1.4 (Sine process). *For $t_1 < t_2$, and compact $E_1$ and $E_2 \subset \mathbb{R}$, the log of the joint probability for the sine processes $S_i(t)$,*

$$
G(t; u_1, \ldots, u_{2r}; v_1, \ldots, v_{2s}) := \log P(\text{all } S_i(t_1) \in E_1^c, \ \text{all } S_i(t_2) \in E_2^c),
$$

*satisfies the third-order nonlinear PDE,*

$$
\begin{aligned}
(1.17) \quad & L_u \frac{(2E_v L_u + (E_v - E_u - 1)L_v)G}{(L_u + L_v)^2 G + \pi^2} \\
& = L_v \frac{(2E_u L_v + (E_u - E_v - 1)L_u)G}{(L_u + L_v)^2 G + \pi^2}.
\end{aligned}
$$

COROLLARY 1.5. *In the case of a single interval, the logarithm of the joint probability for the Sine process,*

$$H(t; x, y) = \log P(S(t_1) \notin [x_1 + x_2, x_1 - x_2], S(t_2) \notin [y_1 + y_2, y_1 - y_2])$$

*satisfies*

$$
\begin{aligned}
(1.18) \quad & \frac{\partial}{\partial x_1} \frac{(2E_y \partial/\partial x_1 + (E_y - E_x - 1)\partial/\partial y_1)H}{(\partial/\partial x_1 + \partial/\partial y_1)^2 H + \pi^2} \\
& = \frac{\partial}{\partial y_1} \frac{(2E_x \partial/\partial y_1 + (E_x - E_y - 1)\partial/\partial x_1)H}{(\partial/\partial x_1 + \partial/\partial y_1)^2 H + \pi^2}.
\end{aligned}
$$

In a recent paper, Tracy and Widom [16] express the joint distribution of the Airy process for several times $t_1, \ldots, t_m$, in terms of an augmented system of auxiliary variables, which satisfy an implicit closed system of nonlinear PDEs. In [17], Tracy and Widom define the Sine process and find an implicit PDE for this process, with methods similar to the one used by them in the Airy process. The quantities involved are entirely different from ours and their methods are functional-theoretical; it remains unclear what the connection is between their results and ours.

The PDEs obtained above provide a very handy tool to compute large time asymptotics for these different processes, with the disadvantage that one usually needs an assumption concerning the interchange of sums and limits; see Section 6 and the Appendix. This is now illustrated for the Airy process, for which we prove the following theorem, assuming some conjecture, mentioned below. This will be discussed in Section 6. A rigorous proof of expansion (1.19) was given later by Widom [18]; his proof was based on the Fredholm determinant expression for the joint distribution.

THEOREM 1.6 (Large time asymptotics for the Airy process). *For large $t = t_2 - t_1$, the joint probability admits the asymptotic series*

$$
\begin{aligned}
(1.19) \quad & P(A(t_1) \leq u, A(t_2) \leq v) \\
& = F_2(u)F_2(v) + \frac{F_2'(u)F_2'(v)}{t^2} + \frac{\Phi(u,v) + \Phi(v,u)}{t^4} + O\left(\frac{1}{t^6}\right),
\end{aligned}
$$

*with the function* $q = q(\alpha)$ *as is the function* (1.6) *and*

$$\Phi(u,v) := F_2(u)F_2(v) \begin{pmatrix} \frac{1}{4}\left(\int_u^\infty q^2\,d\alpha\right)^2 \left(\int_v^\infty q^2\,d\alpha\right)^2 \\ + q^2(u)\left(\frac{1}{4}q^2(v) - \frac{1}{2}\left(\int_v^\infty q^2\,d\alpha\right)^2\right) \\ + \int_v^\infty d\alpha(2(v-\alpha)q^2 + q'^2 - q^4)\int_u^\infty q^2\,d\alpha \end{pmatrix}.$$

*Moreover, the covariance for large* $t = t_2 - t_1$ *behaves as*

(1.20) $\quad E(A(t_2)A(t_1)) - E(A(t_2))E(A(t_1)) = \dfrac{1}{t^2} + \dfrac{c}{t^4} + \cdots,$

*where*

$$c := 2\iint_{\mathbb{R}^2} \Phi(u,v)\,du\,dv.$$

CONJECTURE. The Airy process satisfies the nonexplosion condition for fixed $x$:

(1.21) $\quad \lim_{z\to\infty} P(A(t) \geq x + z \mid A(0) \leq -z) = 0.$

This conjecture will be discussed in Section 6, just before the proof of Theorem 1.6 and in the Appendix.

Finally, in Section 7, we give a rigorous proof of the convergence of the extended Hermite kernel to the Airy and Sine kernels, under the substitutions

$$\mathcal{S}_1 := \left\{ t \mapsto \dfrac{t}{n^{1/3}}, s \mapsto \dfrac{s}{n^{1/3}}, \begin{array}{l} x \mapsto \sqrt{2n+1} + \dfrac{u}{\sqrt{2}n^{1/6}} \\ y \mapsto \sqrt{2n+1} + \dfrac{v}{\sqrt{2}n^{1/6}} \end{array} \right\},$$

$$\mathcal{S}_2 := \left\{ t \mapsto \dfrac{\pi^2 t}{2n}, s \mapsto \dfrac{\pi^2 s}{2n}, x \mapsto \dfrac{\pi u}{\sqrt{2n}}, y \mapsto \dfrac{\pi v}{\sqrt{2n}} \right\}.$$

The precise formula for these kernels will be given later in the beginning of Section 7.

PROPOSITION 1.7. *Under the substitutions* $\mathcal{S}_1$ *and* $\mathcal{S}_2$, *the extended Hermite kernel tends, respectively, to the extended Airy and Sine kernel, when* $n \to \infty$, *uniformly for* $u, v \in$ *compact subsets* $\subset \mathbb{R}$:

$$\lim_{n\to\infty} K_{t,s}^{H,n}(x,y)\,dy_{|\mathcal{S}_1} = K_{t,s}^A(u,v)\,dv,$$

$$\lim_{n\to\infty} K_{t,s}^{H,n}(x,y)\,dy_{|\mathcal{S}_2} = e^{-(\pi^2/2)(t-s)} K_{t,s}^S(u,v)\,dv.$$

**2. The spectrum of coupled random matrices.** Consider a product ensemble $(M_1, M_2) \in \mathcal{H}_n^2 := \mathcal{H}_n \times \mathcal{H}_n$ of $n \times n$ Hermitian matrices, equipped with a Gaussian probability measure,

$$(2.1) \qquad c_n \, dM_1 \, dM_2 \, e^{-\operatorname{Tr}(M_1^2 + M_2^2 - 2cM_1M_2)/2},$$

where $dM_1 \, dM_2$ is the Haar measure on the product $\mathcal{H}_n^2$, with each $dM_i$,

$$(2.2) \qquad dM_1 = \Delta_n^2(x) \prod_1^n dx_i \, dU_1 \quad \text{and} \quad dM_2 = \Delta_n^2(y) \prod_1^n dy_i \, dU_2$$

decomposed into radial and angular parts. In [1], we defined differential operators $\tilde{\mathcal{A}}_k$, $\tilde{\mathcal{B}}_k$ of "weight" $k$, which form a closed Lie algebra, in terms of the coupling constant $c$, appearing in (2.1), and the boundary of the set

$$(2.3) \qquad E = E_1 \times E_2 := \bigcup_{i=1}^r [a_{2i-1}, a_{2i}] \times \bigcup_{i=1}^s [b_{2i-1}, b_{2i}] \subset \mathbb{R}^2.$$

Here we only need the first few ones:

$$(2.4) \qquad \begin{aligned} \tilde{\mathcal{A}}_1 &= \frac{1}{c^2 - 1} \left( \sum_1^{2r} \frac{\partial}{\partial a_j} + c \sum_1^{2s} \frac{\partial}{\partial b_j} \right), \\ \tilde{\mathcal{B}}_1 &= \frac{1}{1 - c^2} \left( c \sum_1^{2r} \frac{\partial}{\partial a_j} + \sum_1^{2s} \frac{\partial}{\partial b_j} \right), \\ \tilde{\mathcal{A}}_2 &= \sum_{j=1}^{2r} a_j \frac{\partial}{\partial a_j} - c \frac{\partial}{\partial c}, \\ \tilde{\mathcal{B}}_2 &= \sum_{j=1}^{2s} b_j \frac{\partial}{\partial b_j} - c \frac{\partial}{\partial c}. \end{aligned}$$

In [1], we proved the following theorem in terms of the Wronskian $\{f, g\}_X = g(Xf) - f(Xg)$, with regard to the first-order differential operator $X$ and based on integrable and Virasoro theory:

THEOREM 2.1. *Given the joint distribution*

$$(2.5) \; P_n(E) := P(\text{all } (M_1\text{-}eigenvalues) \in E_1, \text{all } (M_2\text{-}eigenvalues) \in E_2),$$

*the function* $F_n(c; a_1, \ldots, a_{2r}, b_1, \ldots, b_{2s}) := \log P_n(E)$ *satisfies the nonlinear third-order partial differential equation*

$$(2.6) \qquad \begin{aligned} &\left\{ \tilde{\mathcal{B}}_2 \tilde{\mathcal{A}}_1 F_n, \tilde{\mathcal{B}}_1 \tilde{\mathcal{A}}_1 F_n + \frac{nc}{c^2 - 1} \right\}_{\tilde{\mathcal{A}}_1} \\ &\quad - \left\{ \tilde{\mathcal{A}}_2 \tilde{\mathcal{B}}_1 F_n, \tilde{\mathcal{A}}_1 \tilde{\mathcal{B}}_1 F_n + \frac{nc}{c^2 - 1} \right\}_{\tilde{\mathcal{B}}_1} = 0. \end{aligned}$$

REMARK. Note that both $P_n(E_1 \times E_2)$ and $P_n(E_1^c \times E_2^c)$ satisfy the same equation.

### 3. The joint distribution for the Dyson Brownian motion.

PROOF OF THEOREM 1.1. Using the notation of Section 2 and the change of variables

$$(3.1) \qquad M_i = \frac{B_i}{\sqrt{(1-c^2)/2}},$$

one computes for $t = t_2 - t_1 > 0$, with $e^{-t} = c$,

$$G_n(t; a_1, \ldots, a_{2r}; b_1, \ldots, b_{2s})$$
$$:= \log P(\text{all } \lambda_i(t_1) \in E_1, \text{all } \lambda_i(t_2) \in E_2)$$
$$= \log P(\text{all } (B(t_1)\text{-eigenvalues}) \in E_1, \text{all } (B(t_2)\text{-eigenvalues}) \in E_2)$$
$$= \log \iint_{\substack{\text{all } B_1\text{-eigenvalues } \in E_1 \\ \text{all } B_2\text{-eigenvalues } \in E_2}} Z^{-1} \frac{dB_1 \, dB_2}{(1-c^2)^{n^2/2}} e^{-(1/(1-c^2))\operatorname{Tr}(B_1^2 + B_2^2 - 2cB_1 B_2)}$$
$$= \log \iint_{\substack{\text{all } M_1\text{-eigenvalues } \in E_1/(\sqrt{(1-c^2)/2}) \\ \text{all } M_2\text{-eigenvalues } \in E_2/(\sqrt{(1-c^2)/2})}} Z'^{-1} \, dM_1 \, dM_2$$
$$\times e^{-\operatorname{Tr}(M_1^2 + M_2^2 - 2cM_1 M_2)/2}$$
$$= F_n\bigg(c; \frac{a_1}{\sqrt{(1-c^2)/2}}, \ldots, \frac{a_{2r}}{\sqrt{(1-c^2)/2}};$$
$$\qquad \frac{b_1}{\sqrt{(1-c^2)/2}}, \ldots, \frac{b_{2s}}{\sqrt{(1-c^2)/2}}\bigg),$$

in terms of the function $F_n$ defined in Theorem 2.1. Setting

$$F_n(c; a_1, \ldots, a_{2r}; b_1, \ldots, b_{2s})$$
$$= G_n(t; a_1 \sqrt{(1-c^2)/2}, \ldots, a_{2r} \sqrt{(1-c^2)/2};$$
$$\qquad b_1 \sqrt{(1-c^2)/2}, \ldots, b_{2s} \sqrt{(1-c^2)/2})$$
$$= G_n(t; \tilde{a}_1, \ldots, \tilde{a}_{2r}; \tilde{b}_1, \ldots, \tilde{b}_{2s})$$

in (2.6) leads to the following equation for $G_n := G_n(t; \tilde{a}_1, \ldots, \tilde{a}_{2r}; \tilde{b}_1, \ldots, \tilde{b}_{2s})$, namely,

$$\left\{\tilde{\tilde{\mathcal{B}}}_2 \tilde{\tilde{\mathcal{A}}}_1 G_n, \tilde{\tilde{\mathcal{B}}}_1 \tilde{\tilde{\mathcal{A}}}_1 G_n + \frac{nc}{c^2-1}\right\}_{\tilde{\mathcal{A}}_1} - \left\{\tilde{\tilde{\mathcal{A}}}_2 \tilde{\tilde{\mathcal{B}}}_1 G_n, \tilde{\tilde{\mathcal{A}}}_1 \tilde{\tilde{\mathcal{B}}}_1 G_n + \frac{nc}{c^2-1}\right\}_{\tilde{\mathcal{B}}_1} = 0,$$

where

$$\tilde{\mathcal{A}}_i F_n(c; a_1, \ldots, a_{2r}; b_1, \ldots, b_{2s}) = \tilde{\tilde{\mathcal{A}}}_i G_n(t; \tilde{a}_1, \ldots, \tilde{a}_{2r}; \tilde{b}_1, \ldots, \tilde{b}_{2s}),$$
$$\tilde{\mathcal{B}}_i F_n(c; a_1, \ldots, a_{2r}; b_1, \ldots, b_{2s}) = \tilde{\tilde{\mathcal{B}}}_i G_n(t; \tilde{a}_1, \ldots, \tilde{a}_{2r}; \tilde{b}_1, \ldots, \tilde{b}_{2s}),$$

with the $\tilde{\mathcal{A}}_i, \tilde{\mathcal{B}}_i$ as in (2.4), and

$$
\begin{aligned}
\tilde{\tilde{\mathcal{A}}}_1 &= -\frac{1}{\sqrt{2(1-c^2)}}\left(\sum_1^{2r}\frac{\partial}{\partial \tilde{a}_j} + c\sum_1^{2s}\frac{\partial}{\partial \tilde{b}_j}\right), \\
\tilde{\tilde{\mathcal{B}}}_1 &= \frac{1}{\sqrt{2(1-c^2)}}\left(c\sum_1^{2r}\frac{\partial}{\partial \tilde{a}_j} + \sum_1^{2s}\frac{\partial}{\partial \tilde{b}_j}\right), \\
\tilde{\tilde{\mathcal{A}}}_2 &= \frac{1}{(1-c^2)}\left(\sum_1^{2r}\tilde{a}_j\frac{\partial}{\partial \tilde{a}_j} + c^2\sum_1^{2s}\tilde{b}_j\frac{\partial}{\partial \tilde{b}_j}\right) + \frac{\partial}{\partial t}, \\
\tilde{\tilde{\mathcal{B}}}_2 &= \frac{1}{(1-c^2)}\left(c^2\sum_1^{2r}\tilde{a}_j\frac{\partial}{\partial \tilde{a}_j} + \sum_1^{2s}\tilde{b}_j\frac{\partial}{\partial \tilde{b}_j}\right) + \frac{\partial}{\partial t}.
\end{aligned}
\tag{3.2}
$$

Then clearing the denominators in (3.2) leads to

$$
\{\mathcal{B}_2\mathcal{A}_1 G_n, \mathcal{B}_1\mathcal{A}_1 G_n + 2nc\}_{\mathcal{A}_1} = \{\mathcal{A}_2\mathcal{B}_1 G_n, \mathcal{A}_1\mathcal{B}_1 G_n + 2nc\}_{\mathcal{B}_1},
\tag{3.3}
$$

for the operators (1.10), with $G_n = G_n(t; a_1, \ldots, a_{2r}; b_1, \ldots, b_{2s})$. This establishes Theorem 1.1. □

REMARK. In view of the remark in Section 2, here the expressions

$$
G_n = \begin{cases} \log P(\text{all } \lambda_i(t_1) \in E_1, \text{all } \lambda_i(t_2) \in E_2) \\ \text{or} \\ \log P(\text{all } \lambda_i(t_1) \in E_1^c, \text{all } \lambda_i(t_2) \in E_2^c) \end{cases}
\tag{3.4}
$$

also satisfy the same equation.

## 4. The joint distribution for the Airy process.

PROOF OF THEOREM 1.2. Consider as in (1.13), the disjoint union of intervals

$$
E_1 := \bigcup_{i=1}^r [u_{2i-1}, u_{2i}] \quad \text{and} \quad E_2 := \bigcup_{i=1}^s [v_{2i-1}, v_{2i}] \subseteq \mathbb{R}.
\tag{4.1}
$$

Then, setting $\tau = \tau_2 - \tau_1$,

$$
\begin{aligned}
&H_n(\tau; u_1, \ldots, u_{2r}; v_1, \ldots, v_{2s}) \\
&= \log P\left(\begin{array}{l}\sqrt{2}n^{1/6}(\lambda_n(n^{-1/3}\tau_1) - \sqrt{2n}) \in E_1, \\ \sqrt{2}n^{1/6}(\lambda_n(n^{-1/3}\tau_2) - \sqrt{2n}) \in E_2\end{array}\right) \\
&= \log P(\lambda_n(n^{-1/3}\tau_1) \in \tilde{E}_1, \lambda_n(n^{-1/3}\tau_2) \in \tilde{E}_2) \\
&= G_n(n^{-1/3}\tau; a_1, \ldots, a_{2r}; b_1, \ldots, b_{2s})
\end{aligned}
\tag{4.2}
$$

for the disjoint union of intervals

$$
\tilde{E}_1 := \bigcup_{i=1}^r [a_{2i-1}, a_{2i}] \quad \text{and} \quad \tilde{E}_2 := \bigcup_{i=1}^s [b_{2i-1}, b_{2i}],
\tag{4.3}
$$

with

(4.4) $$a_i = \sqrt{2n} + \frac{u_i}{\sqrt{2}n^{1/6}} \quad \text{and} \quad b_i = \sqrt{2n} + \frac{v_i}{\sqrt{2}n^{1/6}}.$$

The method here is to do asymptotics on (3.3) for $n$ large. In the notation (1.14), define

(4.5) $$L := L_u + L_v, \qquad E := E_u + E_v,$$

with the understanding that $t$ in $E$ now gets replaced by $\tau$. Setting $k := n^{1/6}$ and changing variables,

(4.6)
$$\begin{aligned}
\mathcal{A}_1 G_n(n^{-1/3}\tau; a_1, \ldots, a_{2r}; b_1, \ldots, b_{2s}) \\
= \mathcal{A}_1 H_n(\tau; -k(2k^3 - \sqrt{2}a_1), \ldots; -k(2k^3 - \sqrt{2}b_1), \ldots) \\
= k\sqrt{2}\left(\sum_1^{2r} \frac{\partial}{\partial u_i} + e^{-\tau/k^2} \sum_1^{2r} \frac{\partial}{\partial v_i}\right) H_n(\tau; u_1, \ldots, u_{2r}; v_1, \ldots, v_{2r}) \\
= \mathcal{A}_1|_{\substack{a \to u \\ b \to v}} H_n(\tau; u_1, \ldots, u_{2r}; v_1, \ldots, v_{2r}),
\end{aligned}$$

where the operators $\mathcal{A}_i$ and $\mathcal{B}_i$ are now expressed in $u, v, \tau$-coordinates, using the change of coordinates (4.4) and the chain rule. In these new coordinates, the $\mathcal{A}_i$ and $\mathcal{B}_i$, Taylor expanded in $1/k$, for large $k$, read as follows:

(4.7)
$$\begin{aligned}
\mathcal{A}_1 &= \sqrt{2}k\left(L - \left(\frac{\tau}{k^2} - \frac{\tau^2}{2k^4} + \frac{\tau^3}{6k^6}\right)L_v + O\left(\frac{1}{k^8}\right)\right), \\
\mathcal{B}_1 &= \sqrt{2}k\left(L - \left(\frac{\tau}{k^2} - \frac{\tau^2}{2k^4} + \frac{\tau^3}{6k^6}\right)L_u + O\left(\frac{1}{k^8}\right)\right), \\
\mathcal{A}_2 &= 2k^4\bigg(L - \frac{2\tau}{k^2}L_v + \frac{1}{2k^4}(E - 1 + 4\tau^2 L_v) \\
&\qquad - \frac{\tau}{k^6}\left(E_v - 1 + \frac{4}{3}\tau^2 L_v\right) + O\left(\frac{1}{k^8}\right)\bigg), \\
\mathcal{B}_2 &= 2k^4\bigg(L - \frac{2\tau}{k^2}L_u + \frac{1}{2k^4}(E - 1 + 4\tau^2 L_u) \\
&\qquad - \frac{\tau}{k^6}\left(E_u - 1 + \frac{4}{3}\tau^2 L_u\right) + O\left(\frac{1}{k^8}\right)\bigg).
\end{aligned}$$

Hence, from (4.7),

$$
\begin{aligned}
(4.8)\quad &\frac{1}{2\sqrt{2}k^5}\mathcal{B}_2\mathcal{A}_1 \\
&= L^2 - \frac{\tau}{k^2}(L+L_u)L \\
&\quad + \frac{1}{2k^4}(L(E-2) + \tau^2(4L_u(L+L_v) + LL_v)) \\
&\quad - \frac{\tau}{k^6}\left(L(E_u - 2) + \frac{1}{2}L_v(E+2) + \frac{\tau^2}{6}(8LL_u + 18L_uL_v + LL_v)\right) \\
&\quad + O\!\left(\frac{1}{k^8}\right),
\end{aligned}
$$

$$
\frac{1}{2k^2}\mathcal{B}_1\mathcal{A}_1 = L^2 - \frac{\tau}{k^2}L^2 + \frac{\tau^2}{k^4}\left(\frac{1}{2}L^2 + L_uL_v\right) - \frac{\tau^3}{k^6}\left(\frac{1}{6}L^2 + L_uL_v\right) + O\!\left(\frac{1}{k^8}\right).
$$

As will be shown in Proposition 7.1, using also (7.13), we have for $u = (u_1, \ldots, u_{2r})$ and $v = (v_1, \ldots, v_{2s})$ (also for their derivatives with regard to the endpoints of the intervals)

$$
(4.9)\qquad H_n(\tau; u, v)|_{n=k^6} = G(\tau; u, v) + o(1/k) \qquad \text{for } k \to \infty,
$$

with

$$
G(\tau, u, v) := \log P(A(\tau_1) \in E_1, A(\tau_2) \in E_2),
$$

which will be used below. Proposition 7.1 actually implies that the error in (4.9) has order $k^{-2}\log k$. Equation (3.3), with (4.7) and (4.8) substituted, is a series in $k^2$, for large $k$, but where the three leading coefficients, namely the ones of $k^4$, $k^2$ and $k^0$, vanish:

$$
\begin{aligned}
(4.10)\quad 0 &= \left\{\frac{1}{2\sqrt{2}k^5}\mathcal{B}_2\mathcal{A}_1 H_n, \frac{1}{2k^2}(\mathcal{B}_1\mathcal{A}_1 H_n + 2k^{6e^{-\tau/k^2}})\right\}_{\mathcal{A}_1/(\sqrt{2}k)} \\
&\quad - \left\{\frac{1}{2\sqrt{2}k^5}\mathcal{A}_2\mathcal{B}_1 H_n, \frac{1}{2k^2}(\mathcal{A}_1\mathcal{B}_1 H_n + 2k^6 e^{-\tau/k^2})\right\}_{\mathcal{B}_1/(\sqrt{2}k)} \\
&= \frac{2\tau}{k^2}\bigg[((L_u + L_v)(L_u E_v - L_v E_u) + \tau^2(L_u - L_v)L_uL_v)H_n \\
&\qquad\qquad - \frac{1}{2}\{(L_u^2 - L_v^2)H_n, (L_u+L_v)^2 H_n\}_{L_u+L_v}\bigg] + O\!\left(\frac{1}{k^3}\right) \\
&= \frac{2\tau}{k^2}\bigg[((L_u + L_v)(L_u E_v - L_v E_u) + \tau^2(L_u - L_v)L_uL_v)G \\
&\qquad\qquad - \frac{1}{2}\{(L_u^2 - L_v^2)G, (L_u+L_v)^2 G\}_{L_u+L_v}\bigg] + O\!\left(\frac{1}{k^3}\right),
\end{aligned}
$$

using (4.9). In this calculation, we used the linearity of the Wronskian $\{X, Y\}_Z$ in the three arguments and the following commutation relations:

$$
(4.11)\quad [L_u, E_u] = L_u, \qquad [L_u, E_v] = [L_u, L_v] = [L_u, \tau] = 0, \qquad [E_u, \tau] = \tau,
$$

including their dual relations by $u \leftrightarrow v$; also $\{L^2 G, 1\}_{L_u - L_v} = \{L(L_u - L_v)G, 1\}_L$. It is also useful to note in (4.10), that the two Wronskians in the first expression are dual to each other by $u \leftrightarrow v$. The point of the computation is to preserve the Wronskian structure up to the end. This proves Theorem 1.2. □

PROOF OF COROLLARY 1.3.  Equation (1.15) for the probability

$$G(\tau; u, v) := \log P(A(\tau_1) \leq u, A(\tau_2) \leq v), \qquad \tau = \tau_2 - \tau_1,$$

takes on the explicit form

(4.12)
$$\begin{aligned}
\tau \frac{\partial}{\partial \tau}\left(\frac{\partial^2}{\partial u^2} - \frac{\partial^2}{\partial v^2}\right)G &= \frac{\partial^3 G}{\partial u^2 \, \partial v}\left(2\frac{\partial^2 G}{\partial v^2} + \frac{\partial^2 G}{\partial u \, \partial v} - \frac{\partial^2 G}{\partial u^2} + u - v - \tau^2\right) \\
&\quad - \frac{\partial^3 G}{\partial v^2 \, \partial u}\left(2\frac{\partial^2 G}{\partial u^2} + \frac{\partial^2 G}{\partial u \, \partial v} - \frac{\partial^2 G}{\partial v^2} - u + v - \tau^2\right) \\
&\quad + \left(\frac{\partial^3 G}{\partial u^3}\frac{\partial}{\partial v} - \frac{\partial^3 G}{\partial v^3}\frac{\partial}{\partial u}\right)\left(\frac{\partial}{\partial u} + \frac{\partial}{\partial v}\right)G.
\end{aligned}$$

This equation enjoys an obvious $u \leftrightarrow v$ duality. Finally the change of variables in the statement of Corollary 1.3 leads to (1.16). □

## 5. The joint distribution for the Sine process.

PROOF OF THEOREM 1.4.  Consider as in (1.13), the disjoint union of (in this case) compact intervals $E_1 := \bigcup_{i=1}^r [u_{2i-1}, u_{2i}]$ and $E_2 := \bigcup_{i=1}^s [v_{2i-1}, v_{2i}] \subseteq \mathbb{R}$. Then, again setting $\tau = \tau_2 - \tau_1$,

(5.1)
$$\begin{aligned}
&H_n(\tau; u_1, \ldots, u_{2r}; v_1, \ldots, v_{2s}) \\
&= \log P\left(\text{all } 2\sqrt{n}\lambda_{(n/2)+i}\left(\frac{\tau_1}{n}\right) \in E_1^c, \text{ all } 2\sqrt{n}\lambda_{(n/2)+i}\left(\frac{\tau_2}{n}\right) \in E_2^c\right) \\
&= \log P\left(\text{all } \lambda_{(n/2)+i}\left(\frac{\tau_1}{n}\right) \in \tilde{E}_1^c, \text{all } \lambda_{(n/2)+i}\left(\frac{\tau_2}{n}\right) \in \tilde{E}_2^c\right) \\
&= G_n\left(\frac{\tau}{n}; a_1, \ldots, a_{2r}; b_1, \ldots, b_{2s}\right)
\end{aligned}$$

for the disjoint union of intervals

(5.2) $$\tilde{E}_1 := \bigcup_{i=1}^r [a_{2i-1}, a_{2i}] \quad \text{and} \quad \tilde{E}_2 := \bigcup_{i=1}^s [b_{2i-1}, b_{2i}],$$

with

(5.3) $$a_i := \frac{u_i}{2\sqrt{n}}, \qquad b_i := \frac{v_i}{2\sqrt{n}}.$$

Note here $G_n$ refers to the second formula in (3.4). Setting $k := n^{1/2}$, using the change of variables (5.3) and the chain rule,

$$\mathcal{A}_1 G_n(\tau/k^2; a_1, \ldots, a_{2r}; b_1, \ldots, b_{2s})$$
$$= \mathcal{A}_1 H_n(\tau; 2ka_1, \ldots; 2kb_1, \ldots)$$
$$= 2k\left(\sum_1^{2r} \frac{\partial}{\partial u_i} + e^{-\tau/k^2} \sum_1^{2r} \frac{\partial}{\partial v_i}\right) H_n(\tau; u_1, \ldots, u_{2r}; v_1, \ldots, v_{2r})$$
$$= \mathcal{A}_{1|\substack{a \to u \\ b \to v}} H_n(\tau; u_1, \ldots, u_{2r}; v_1, \ldots, v_{2r}).$$

In these new $u, v, \tau$-coordinates, the operators $\mathcal{A}_i$ and $\mathcal{B}_i$, Taylor expanded in $1/k$ for large $k$, read as follows:

(5.4)
$$\mathcal{A}_1 = 2k\left(L - \frac{\tau}{k^2}L_v + O\left(\frac{1}{k^4}\right)\right),$$
$$\mathcal{B}_1 = 2k\left(L - \frac{\tau}{k^2}L_u + O\left(\frac{1}{k^4}\right)\right),$$
$$\mathcal{A}_2 = E - 1 - \frac{2\tau}{k^2}(E_v - 1) + O\left(\frac{1}{k^4}\right),$$
$$\mathcal{B}_2 = E - 1 - \frac{2\tau}{k^2}(E_u - 1) + O\left(\frac{1}{k^4}\right).$$

Moreover, as will be shown in Proposition 7.3 (Section 7.3) we have the following asymptotic estimate (and similarly for its derivatives with regard to the endpoints of the intervals):

(5.5) $$H_n(\tau; u, v)|_{n=k^2} = \tilde{G}(\tau; u, v) + o(1) \qquad \text{for } k \to \infty,$$

with

(5.6) $$\tilde{G}(\tau, u, v) := \log P\left(\text{all } \sqrt{2}\pi S_i\left(\frac{2\tau_1}{\pi^2}\right) \in E_1^c, \text{all } \sqrt{2}\pi S_i\left(\frac{2\tau_2}{\pi^2}\right) \in E_2^c\right).$$

Using the expansions (5.4) of the $\mathcal{A}_i$ and $\mathcal{B}_i$ and later the commutation relations (4.11), yields the following for the Wronskian:

$$\frac{1}{(2k)^4}\{\mathcal{B}_2\mathcal{A}_1 H_n, \mathcal{B}_1\mathcal{A}_1 H_n + 2k^2 e^{-\tau/k^2}\}_{\mathcal{A}_1}$$
$$= \frac{1}{(2k)^4}\left\{2k\left(\left(E - 1 - \frac{2\tau}{k^2}(E_u - 1)\right)\left(L - \frac{\tau}{k^2}L_v\right) + O\left(\frac{1}{k^4}\right)\right)H_n,\right.$$
$$(2k)^2\left(\left(L - \frac{\tau}{k^2}L_v\right)\left(L - \frac{\tau}{k^2}L_u\right) + O\left(\frac{1}{k^4}\right)\right)H_n$$
$$\left. + \frac{(2k)^2}{2}\left(1 - \frac{\tau}{k^2} + O\left(\frac{1}{k^4}\right)\right)\right\}_{2k(L-(\tau/k^2)L_v + O(1/k^4))}$$
$$= \left\{\left((E-1)L - \frac{2\tau}{k^2}(E_u - 1)L - \frac{1}{k^2}(E-1)\tau L_v + O\left(\frac{1}{k^4}\right)\right)H_n,\right.$$

$$L^2 H_n + \frac{1}{2} - \frac{\tau}{k^2}\left(L^2 H_n + \frac{1}{2}\right) + O\left(\frac{1}{k^4}\right)\bigg\}_{L-(\tau/k^2)L_v}$$

$$= \bigg\{\left((E-1)L - \frac{\tau}{k^2}(2(E_u-1)L + (E+1)L_v) + O\left(\frac{1}{k^4}\right)\right)H_n,$$

$$L^2 H_n + \frac{1}{2} - \frac{\tau}{k^2}\left(L^2 H_n + \frac{1}{2}\right) + O\left(\frac{1}{k^4}\right)\bigg\}_{L-(\tau/k^2)L_v}$$

$$= \left\{(E-1)LH_n, L^2 H_n + \frac{1}{2}\right\}_L - \frac{\tau}{k^2}\left\{(E-1)LH_n, L^2 H_n + \frac{1}{2}\right\}_L$$

$$- \frac{\tau}{k^2}\left\{(E-1)LH_n, L^2 H_n + \frac{1}{2}\right\}_{L_v}$$

$$- \frac{\tau}{k^2}\left\{(2(E_u-1)L + (E+1)L_v)H_n, L^2 H_n + \frac{1}{2}\right\}_L + O\left(\frac{1}{k^4}\right).$$

Hence, subtracting the previous formula from its dual and using (5.5),

$$0 = \frac{k^2}{\tau}\frac{1}{(2k)^4}(\{\mathcal{B}_2\mathcal{A}_1 H_n, \mathcal{B}_1\mathcal{A}_1 H_n + 2k^2 e^{-\tau/k^2}\}_{\mathcal{A}_1}$$

$$- \{\mathcal{A}_2\mathcal{B}_1 H_n, \mathcal{A}_1\mathcal{B}_1 H_n + 2k^2 e^{-\tau/k^2}\}_{\mathcal{B}_1})$$

$$= \left\{(E-1)LH_n, L^2 H_n + \frac{1}{2}\right\}_{L_u - L_v}$$

$$- \left\{(2(E_u-1)L + (E+1)L_v)H_n, LH_n + \frac{1}{2}\right\}_L$$

$$+ \left\{(2(E_v-1)L + (E+1)L_u)H_n, LH_n + \frac{1}{2}\right\}_L + O\left(\frac{1}{k^2}\right)$$

$$= \left\{(E-1)L\tilde{G}, L^2\tilde{G} + \frac{1}{2}\right\}_{L_u - L_v}$$

$$- \left\{(2(E_u-1)L + (E+1)L_v)\tilde{G}, L\tilde{G} + \frac{1}{2}\right\}_L$$

$$+ \left\{(2(E_v-1)L + (E+1)L_u)\tilde{G}, L\tilde{G} + \frac{1}{2}\right\}_L + o(1)$$

$$= \left\{(E-1)L\tilde{G}, L^2\tilde{G} + \frac{1}{2}\right\}_{L_u - L_v}$$

$$+ \left\{(2(E_v - E_u)L + (E+1)(L_u - L_v))\tilde{G}, L^2\tilde{G} + \frac{1}{2}\right\}_L + o(1).$$

Upon division by $(L^2\tilde{G} + \frac{1}{2})^2$, one finds

$$0 = (L_u - L_v)\left(\frac{(E-1)L\tilde{G}}{L^2\tilde{G}+1/2}\right)$$
$$+ (L_u + L_v)\left(\frac{(2(E_v - E_u)L + (E+1)(L_u - L_v))\tilde{G}}{L^2\tilde{G}+1/2}\right)$$

$$= L_u \frac{((E-1)L + 2(E_v - E_u)L + (E+1)(L_u - L_v))\tilde{G}}{L^2\tilde{G}+1/2}$$

$$+ L_v \frac{(-(E-1)L + 2(E_v - E_u)L + (E+1)(L_u - L_v))\tilde{G}}{L^2\tilde{G}+1/2}$$

$$= 2L_u \frac{(2E_v L_u + (E_v - E_u - 1)L_v)\tilde{G}}{L^2\tilde{G}+1/2}$$

$$- 2L_v \frac{(2E_u L_v + (E_u - E_v - 1)L_u)\tilde{G}}{L^2\tilde{G}+1/2}.$$

In view of (5.6), the function

$$G(\tau, u, v) = \log P(\text{all } S_i(\tau_1) \in E_1^c, \text{all } S_i(\tau_2) \in E_2^c)$$
$$= \tilde{G}\left(\frac{\pi^2\tau}{2}, \sqrt{2}\pi u, \sqrt{2}\pi v\right)$$

satisfies (1.17) of Theorem 1.4. □

PROOF OF COROLLARY 1.5. Setting

$$u_1 = x_1 + x_2, \qquad u_2 = x_1 - x_2,$$
$$v_1 = y_1 + y_2, \qquad v_2 = y_1 - y_2,$$

the function

(5.7) $\quad H(t; x_1, x_2; y_1, y_2) := G(t; x_1 + x_2, x_1 - x_2; y_1 + y_2, y_1 - y_2)$

satisfies (1.18), ending the proof of Corollary 1.5. □

**6. Large time asymptotics for the Airy process.** This section aims at proving Theorem 1.6, for which we need the following lemma:

LEMMA 6.1. *The following ratio of probabilities admits the asymptotic expansion for large $t > 0$ in terms of functions $f_i(u, v)$, symmetric in $u$ and*

$v$:

$$\text{(6.1)} \qquad \frac{P(A(0) \leq u, A(t) \leq v)}{P(A(0) \leq u)P(A(t) \leq v)} = 1 + \sum_{i \geq 1} \frac{f_i(u,v)}{t^i},$$

*from which it follows that*

$$\lim_{t \to \infty} P(A(0) \leq u, A(t) \leq v) = P(A(0) \leq u)P(A(t) \leq v) = F_2(u)F_2(v);$$

*this means the Airy process decouples at $\infty$.*

PROOF. This will be done in Section 7.5, using the extended Airy kernel. Note, since the probabilities in (6.1) are symmetric in $u$ and $v$, the coefficients $f_i$ are symmetric as well. The last equality in the formula above follows from stationarity. $\square$

CONJECTURE. *The coefficients $f_i(u,v)$ have the property*

$$\text{(6.2)} \qquad \lim_{u \to \infty} f_i(u,v) = 0 \qquad \text{for fixed } v \in \mathbb{R}$$

*and*

$$\text{(6.3)} \qquad \lim_{z \to \infty} f_i(-z, z+x) = 0 \qquad \text{for fixed } x \in \mathbb{R}.$$

The justification for this plausible conjecture will now follow: First, considering the following conditional probability:

$$P(A(t) \leq v \mid A(0) \leq u) = \frac{P(A(0) \leq u, A(t) \leq v)}{P(A(0) \leq u)}$$

$$= F_2(v)\left(1 + \sum_{i \geq 1} \frac{f_i(u,v)}{t^i}\right),$$

and letting $v \to \infty$, we have automatically

$$1 = \lim_{v \to \infty} P(A(t) \leq v \mid A(0) \leq u) = \lim_{v \to \infty}\left[F_2(v)\left(1 + \sum_{i \geq 1} \frac{f_i(u,v)}{t^i}\right)\right]$$

$$= 1 + \lim_{v \to \infty} \sum_{i \geq 1} \frac{f_i(u,v)}{t^i},$$

which would imply, assuming the interchange of the limit and the summation is valid,

$$\text{(6.4)} \qquad \lim_{v \to \infty} f_i(u,v) = 0$$

and, by symmetry

$$\lim_{u \to \infty} f_i(u,v) = 0.$$

To deal with (6.3), we assume the following *nonexplosion* condition, whose plausibility is discussed in the Appendix: for any fixed $t > 0$, $x \in \mathbb{R}$, the conditional probability satisfies

(6.5) $$\lim_{z \to \infty} P(A(t) \geq x + z \mid A(0) \leq -z) = 0.$$

Hence, the conditional probability satisfies, upon setting

$$v = z + x, \qquad u = -z,$$

and using $\lim_{z \to \infty} F_2(z + x) = 1$, the following:

$$1 = \lim_{z \to \infty} P(A(t) \leq z + x \mid A(0) \leq -z) = 1 + \lim_{z \to \infty} \sum_{i \geq 1} \frac{f_i(-z, z + x)}{t^i},$$

which, assuming the validity of the same interchange, implies that

$$\lim_{z \to \infty} f_i(-z, z + x) = 0 \qquad \text{for all } i \geq 1.$$

PROOF OF THEOREM 1.6. Putting the log of the expansion (6.1)

(6.6) $$\begin{aligned} G(t; u, v) &= \log P(A(0) \leq u, A(t) < v) \\ &= \log F_2(u) + \log F_2(v) + \sum_{i \geq 1} \frac{h_i(u, v)}{t^i} \\ &= \log F_2(u) + \log F_2(v) \\ &\quad + \frac{f_1(u, v)}{t} + \frac{f_2(u, v) - f_1^2(u, v)/2}{t^2} + \cdots, \end{aligned}$$

in (4.12), leads to:

(i) a leading term of order $t$, given by

(6.7) $$\mathcal{L} h_1 = 0,$$

where

(6.8) $$\mathcal{L} := \left( \frac{\partial}{\partial u} - \frac{\partial}{\partial v} \right) \frac{\partial^2}{\partial u \, \partial v}.$$

The most general solution to (6.7) is given by

$$h_1(u, v) = r_1(u) + r_3(v) + r_2(u + v),$$

with arbitrary functions $r_1, r_2, r_3$. Hence,

$$P(A(0) \leq u, A(t) \leq v) = F_2(u) F_2(v) \left( 1 + \frac{h_1(u, v)}{t} + \cdots \right)$$

with $h_1(u, v) = f_1(u, v)$ as in (6.1). Applying (6.2),

$$r_1(u) + r_3(\infty) + r_2(\infty) = 0 \qquad \text{for all } u \in \mathbb{R},$$

implying
$$r_1(u) = \text{constant} = r_1(\infty),$$
and similarly
$$r_3(u) = \text{constant} = r_3(\infty).$$

Therefore, without loss of generality, we may absorb the constants $r_1(\infty)$ and $r_3(\infty)$ in the definition of $r_2(u+v)$. Hence, from (6.6),
$$f_1(u,v) = h_1(u,v) = r_2(u+v)$$
using (6.3),
$$0 = \lim_{z \to \infty} f_1(-z, z+x) = r_2(x)$$
implying that the $h_1(u,v)$-term in the series (6.6) vanishes.

(ii) One computes that the term $h_2(u,v)$ in the expansion (6.6) of $G(t;u,v)$ satisfies
$$(6.9) \qquad \mathcal{L} h_2 = \frac{\partial^3 g}{\partial u^3} \frac{\partial^2 g}{\partial v^2} - \frac{\partial^3 g}{\partial v^3} \frac{\partial^2 g}{\partial u^2} \qquad \text{with } g(u) := \log F_2(u).$$

This is the term of order $t^0$, by putting the series (6.6) in (4.12). The most general solution to (6.9) is
$$h_2(u,v) = g'(u)g'(v) + r_1(u) + r_3(v) + r_2(u+v).$$
Then
$$P(A(0) \le u, A(t) \le v) = e^{G(t;u,v)}$$
$$= F_2(u) F_2(v) \exp \sum_{i \ge 2} \frac{h_i(u,v)}{t^i}$$
$$= F_2(u) F_2(v) \left(1 + \frac{h_2(u,v)}{t^2} + \cdots\right).$$

In view of the explicit formula for the distribution $F_2$ and the behavior (1.6) of $q(\alpha)$ for $\alpha \nearrow \infty$, we have that
$$\lim_{u \to \infty} g'(u) = \lim_{u \to \infty} (\log F_2(u))'$$
$$= \lim_{u \to \infty} \int_u^\infty q^2(\alpha)\, d\alpha = 0.$$
Hence
$$0 = \lim_{u \to \infty} f_2(u,v) = \lim_{u \to \infty} h_2(u,v) = r_1(\infty) + r_3(v) + r_2(\infty),$$

showing $r_3$ and similarly $r_1$ are constants. Therefore, by absorbing $r_1(\infty)$ and $r_3(\infty)$ into $r_2(u+v)$, we have

$$f_2(u,v) = h_2(u,v) = g'(u)g'(v) + r_2(u+v).$$

Again, by the behavior of $q(x)$ at $+\infty$ and $-\infty$, we have for large $z > 0$,

$$g'(-z)g'(z+x) = \int_{-z}^{\infty} q^2(\alpha)\,d\alpha \int_{z+x}^{\infty} q^2(\alpha)\,d\alpha \leq c z^{3/2} e^{-2z/3}.$$

Hence

$$0 = \lim_{z \to \infty} f_2(-z, z+x) = r_2(x)$$

and so

$$f_2(u,v) = h_2(u,v) = g'(u)g'(v),$$

yielding the $1/t^2$ term in the series (6.6).

(iii) Next, setting

(6.10)
$$\begin{aligned}G(t;u,v) &= \log P(A(0) \leq u, A(t) \leq v) \\ &= g(u) + g(v) + \frac{g'(u)g'(v)}{t^2} + \frac{h_3(u,v)}{t^3} + \cdots\end{aligned}$$

in (4.12), we find for the $t^{-1}$ term:

$$\mathcal{L}h_3 = 0.$$

As in (6.7), its most general solution is given by

$$h_3(u,v) = r_1(u) + r_3(v) + r_2(u+v).$$

By exponentiation of (6.6), we find

$$\begin{aligned}P(A(0) &\leq u, A(t) \leq v) \\ &= F_2(u)F_2(v)\left(1 + \frac{g'(u)g'(v)}{t^2} + \frac{r_1(u) + r_3(v) + r_2(u+v)}{t^3} + \cdots\right).\end{aligned}$$

The precise same arguments lead to $h_3(u,v) = 0$.

(iv) So, at the next stage, we have, remembering $g(u) = \log F_2(u)$,

(6.11) $$G(t;u,v) = g(u) + g(v) + \frac{g'(u)g'(v)}{t^2} + \frac{h_4(u,v)}{t^4} + \cdots$$

with

(6.12) $$f_4(u,v) = h_4(u,v) + \tfrac{1}{2}h_2^2(u,v) = h_4(u,v) + \tfrac{1}{2}g'(u)^2 g'(v)^2.$$

Setting the series (6.11) in (4.12), we find for the $t^{-2}$ term:

$$\begin{aligned}
\mathcal{L}h_4 &= 2\left(\frac{\partial^3 g}{\partial u^3}\left(\frac{\partial^2 g}{\partial v^2}\right)^2 - \frac{\partial^3 g}{\partial v^3}\left(\frac{\partial^2 g}{\partial u^2}\right)^2\right) + \frac{\partial^3 g}{\partial u^3}\frac{\partial^3 g}{\partial v^3}\left(\frac{\partial g}{\partial u} - \frac{\partial g}{\partial v}\right) \\
&\quad + \frac{1}{2}\left(\frac{\partial^4 g}{\partial u^4}\frac{\partial}{\partial v}\left(\frac{\partial g}{\partial v}\right)^2 - \frac{\partial^4 g}{\partial v^4}\frac{\partial}{\partial u}\left(\frac{\partial g}{\partial u}\right)^2\right) \\
&\quad + \left(\frac{\partial^3 g}{\partial u^3}\frac{\partial^2 g}{\partial v^2} + \frac{\partial^3 g}{\partial v^3}\frac{\partial^2 g}{\partial u^2}\right)(u-v) + 2\left(\frac{\partial^3 g}{\partial u^3}\frac{\partial g}{\partial v} - \frac{\partial^3 g}{\partial v^3}\frac{\partial g}{\partial u}\right) \\
&= 2(2q(u)q'(u)(q(v)q'(v)+1) \\
&\quad - q(u)q''(u)q^2(v) - (q'(u))^2 q^2(v))\int_v^\infty q^2 \\
&\quad + 2q(u)(q(u)q'(v)q''(v) + q'(u)q(v)q''(v) - 2q(u)q^3(v)q'(v)) \\
&\quad - \text{ same with } u \leftrightarrow v.
\end{aligned}$$
(6.13)

The latter is an expression in $q(u)$, $q(v)$ and its derivatives and in $\int_u^\infty q^2(\alpha)\,d\alpha$ and $\int_v^\infty q^2(\alpha)\,d\alpha$. It is obtained by substituting in the previous expression

$$g(u) = \int_u^\infty (u-\alpha)q^2(\alpha)\,d\alpha$$

and the Painlevé II differential equation for $q(u)$,

$$uq(u) = q''(u) - 2q(u)^3,$$

in order to eliminate the explicit appearance of $u$ and $v$. Now introduce

$$g(u) = \int_u^\infty (u-\alpha)q^2(\alpha)\,d\alpha,$$

$$g_1(u) = \int_u^\infty (u-\alpha)q'^2(\alpha)\,d\alpha,$$

$$g_2(u) = \int_u^\infty (u-\alpha)q^4(\alpha)\,d\alpha.$$

Note $g'(u) = \int_u^\infty q^2(\alpha)\,d\alpha$, $g''(u) = -q^2(u)$ and $g_1'(u) = \int_u^\infty q'^2(\alpha)\,d\alpha$, $g_2'(u) = \int_u^\infty q^4(\alpha)\,d\alpha$. The most general solution to (6.13) is given, modulo the null-space of $\mathcal{L}$, by

$$\begin{aligned}
h_4(u,v) &= \frac{1}{2}(g''(u)g'(v)^2 + g''(v)g'(u)^2 + g''(u)g''(v)) \\
&\quad + g'(u)(2g(v) + g_1'(v) - g_2'(v)) \\
&\quad + g'(v)(2g(u) + g_1'(u) - g_2'(u)) \\
&= q^2(u)\left(\frac{q^2(v)}{4} - \frac{1}{2}\left(\int_v^\infty q^2(\alpha)\,d\alpha\right)^2\right) \\
&\quad + \int_u^\infty q^2(\alpha)\,d\alpha \int_v^\infty (2(v-\alpha)q^2(\alpha) + q'^2(\alpha) - q^4(\alpha))\,d\alpha \\
&\quad + \text{same with } u \leftrightarrow v.
\end{aligned}$$
(6.14)

This form, together with (6.12), implies for the function $f_4(u,v)$:

$$f_4(u,v) = h_4(u,v) + \tfrac{1}{2}g'(u)^2 g'(v)^2 + r_1(u) + r_3(v) + r_2(u+v)$$
$$= \sum_i a_i(u)b_i(v) + r_1(u) + r_3(v) + r_2(u+v).$$

Using the asymptotics of $q(u)$, one finds

$$a_i(u), b_i(u) \leq \begin{cases} ce^{-u}, & u \to \infty, \\ c|u|^3, & u \to -\infty, \end{cases}$$

and so, by the same argument,

$$r_1(u) = r_2(u) = r_3(u) = 0.$$

Therefore, we have

$$f_4(u,v) = h_4(u,v) + \tfrac{1}{2}g'(u)^2 g'(v)^2,$$

with $h_4(u,v)$ as in (6.14), thus yielding (1.19).

Finally, to prove (1.20), we compute from (1.19), after integration by parts, taking into account the boundary terms, using (1.6),

$$E(A(0)A(t)) = \iint_{\mathbb{R}^2} uv \frac{\partial^2}{\partial u\, \partial v} P(A(0) \leq u, A(t) \leq v)\, du\, dv$$
$$= \int_{-\infty}^{\infty} u F_2'(u)\, du \int_{-\infty}^{\infty} v F_2'(v)\, dv$$
$$+ \frac{1}{t^2} \int_{-\infty}^{\infty} F_2'(u)\, du \int_{-\infty}^{\infty} F_2'(v)\, dv$$
$$+ \frac{1}{t^4} \iint_{\mathbb{R}^2} (\Phi(u,v) + \Phi(v,u))\, du\, dv$$
$$+ O\left(\frac{1}{t^6}\right)$$
$$= (E(A(0)))^2 + \frac{1}{t^2} + \frac{c}{t^4} + O\left(\frac{1}{t^6}\right),$$

where

$$c := \iint_{\mathbb{R}^2} (\Phi(u,v) + \Phi(v,u))\, du\, dv = 2 \iint_{\mathbb{R}^2} \Phi(u,v)\, du\, dv,$$

thus ending the proof of Theorem 1.6. $\square$

**7. The extended kernels.** The joint probabilities for the Dyson, Airy and Sine processes can also be expressed in terms of the Fredholm determinant of matrix kernels, the so-called extended Hermite, Airy and Sine kernels

(considered in [6], [11] and especially in [12] and [8]), defined for subsets $E_i \subset \mathbb{R}$,

$$(7.1) \qquad \hat{K}_{t_i t_j}(x,y) := \chi_{E_i^c}(x) K_{t_i t_j}(x,y) \chi_{E_j^c}(y)$$

with $K_{t_i t_j}$ being one of the following kernels:

$$K_{t_i t_j}^{H,n}(x,y) := \begin{cases} \sum_{k=1}^{\infty} e^{-k(t_i-t_j)} \varphi_{n-k}(x) \varphi_{n-k}(y), & \text{if } t_i \geq t_j, \\ -\sum_{k=-\infty}^{0} e^{k(t_j-t_i)} \varphi_{n-k}(x) \varphi_{n-k}(y), & \text{if } t_i < t_j, \end{cases}$$

$$K_{t_i t_j}^{A}(x,y) := \begin{cases} \int_0^{\infty} e^{-z(t_i-t_j)} \mathrm{Ai}(x+z) \mathrm{Ai}(y+z)\,dz, & \text{if } t_i \geq t_j, \\ -\int_{-\infty}^{0} e^{z(t_j-t_i)} \mathrm{Ai}(x+z) \mathrm{Ai}(y+z)\,dz, & \text{if } t_i < t_j, \end{cases}$$

$$K_{t_i t_j}^{S}(x,y) := \begin{cases} \dfrac{1}{\pi} \int_0^{\pi} e^{z^2(t_i-t_j)/2} \cos z(x-y)\,dz, & \text{if } t_i \geq t_j, \\ -\dfrac{1}{\pi} \int_{\pi}^{\infty} e^{-z^2(t_j-t_i)/2} \cos z(x-y)\,dz, & \text{if } t_i < t_j, \end{cases}$$

where

$$(7.2)\ \varphi_k(x) = \begin{cases} e^{-x^2/2} p_k(x), & \text{for } k \geq 0, \text{ with } p_k(x) = \dfrac{H_k(x)}{2^{k/2}\sqrt{k!}\pi^{1/4}}, \\ 0, & \text{for } k < 0; \end{cases}$$

$p_k(x)$ are the normalized Hermite polynomials, and $\mathrm{Ai}(x)$ is the Airy function. Now we make a few comments about these kernels.

7.1. *The Fredholm determinant of extended kernels.* Letting $x(t)$ denote either the largest eigenvalue $\lambda_n(t)$ in the Dyson process, or the Airy process $A(t)$, or the collection of $S_i(t)$'s in the Sine process, the probability is now defined by (drop the superscripts in $\hat{K}^{H,n}, \hat{K}^A, \hat{K}^S$)

$$\begin{aligned}
P(x(t_i) &\in E_i, 1 \leq i \leq m) \\
&= \det(I - z(\hat{K}_{t_i t_j})_{1 \leq i,j \leq m})|_{z=1} \\
&= 1 + \sum_{N=1}^{\infty} (-z)^N \sum_{\substack{0 \leq r_i \leq N \\ \sum_1^m r_i = N}} \int_{\mathcal{R}} \prod_1^{r_1} d\alpha_i^{(1)} \cdots \prod_1^{r_m} d\alpha_i^{(m)} \\
&\quad \times \det\left((\hat{K}_{t_k t_\ell}(\alpha_i^{(k)}, \alpha_j^{(\ell)}))_{\substack{1 \leq i \leq r_k \\ 1 \leq j \leq r_\ell}}\right)_{1 \leq k,\ell \leq m}\bigg|_{z=1},
\end{aligned}$$

where the $N$-fold integral above is taken over the range
$$\mathcal{R} = \left\{\begin{array}{c} -\infty < \alpha_1^{(1)} \leq \cdots \leq \alpha_{r_1}^{(1)} < \infty \\ \vdots \\ -\infty < \alpha_1^{(m)} \leq \cdots \leq \alpha_{r_m}^{(m)} < \infty \end{array}\right\},$$

with integrand equal to the determinant of an $N \times N$ matrix, with blocks given by the $r_k \times r_\ell$ matrices $(\hat{K}_{t_k t_\ell}(\alpha_i^{(k)}, \alpha_j^{(\ell)}))_{\substack{1 \leq i \leq r_k \\ 1 \leq j \leq r_\ell}}$. In particular, for $m = 2$, we have

$$P(x(t_1) \in E_1, x(t_2) \in E_2)$$
$$= 1 + \sum_{N=1}^{\infty} (-z)^N \sum_{\substack{0 \leq r,s \leq N \\ r+s=N}} \int_{\left\{\begin{array}{c} -\infty < \alpha_1 \leq \cdots \leq \alpha_r < \infty \\ -\infty < \beta_1 \leq \cdots \leq \beta_s < \infty \end{array}\right\}} \prod_1^r d\alpha_i \prod_1^s d\beta_i$$
$$\times \det \begin{pmatrix} (\hat{K}_{t_1 t_1}(\alpha_i, \alpha_j))_{1 \leq i,j \leq r} & (\hat{K}_{t_1 t_2}(\alpha_i, \beta_j))_{\substack{1 \leq i \leq r \\ 1 \leq j \leq s}} \\ (\hat{K}_{t_2 t_1}(\beta_i, \alpha_j))_{\substack{1 \leq i \leq s \\ 1 \leq j \leq r}} & (\hat{K}_{t_2 t_2}(\beta_i, \beta_j))_{1 \leq i,j \leq s} \end{pmatrix}\bigg|_{z=1}.$$

**7.2. The extended Hermite kernel tends to the extended Airy kernel.** Given the substitution

(7.3) $$\mathcal{S}_1 := \left\{ t \mapsto \frac{t}{n^{1/3}}, s \mapsto \frac{s}{n^{1/3}}, \begin{array}{c} x \mapsto \sqrt{2n+1} + \dfrac{u}{\sqrt{2}n^{1/6}} \\ y \mapsto \sqrt{2n+1} + \dfrac{v}{\sqrt{2}n^{1/6}} \end{array} \right\},$$

we have:

PROPOSITION 7.1. *The extended Hermite kernel tends to the extended Airy kernel, when $n \to \infty$, uniformly for $u, v \in$ compact subsets $\subset \mathbb{R}$:*

$$|K_{t,s}^{H,n}(x,y) \, dy_{|\mathcal{S}_1} - K_{t,s}^A(u,v) \, dv| = O\left(\frac{\log n}{n^{1/3}}\right).$$

Before proving this proposition, we need the following estimate:

LEMMA 7.2. *For large $n > 0$ and $-M_0 n^{1/3} \log n \leq k \leq M_0 n^{1/3} \log n$, with fixed $M_0 > 0$, we have*

(7.4) $$\varphi_{n-k}\left(\sqrt{2n+1} + \frac{u}{\sqrt{2}n^{1/6}}\right) = \frac{2^{1/4}}{n^{1/12}} \text{Ai}\left(u + \frac{k}{n^{1/3}}\right)(1 + E_u(k,n)),$$

*with the following uniform bound in $u \in$ compact subsets $\mathbb{R}$:*

$$|E_u(k,n)| \leq O\left(\frac{\log n}{n^{2/3}}\right).$$

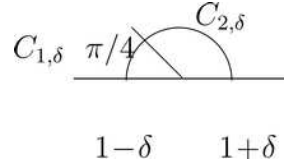

Fig. 1.

PROOF. Here one needs the asymptotics for the Hermite polynomials when $z \in C_{2,\delta}$, as in Figure 1; it is given by (see [3]):

$$\varphi_n(z\sqrt{2n}) = p_n(z\sqrt{2n})e^{-nz^2}$$
$$(7.5) \qquad = \frac{1+O(1/n)}{(2n)^{1/4}}\left\{\left(\frac{z+1}{z-1}f_n(z)\right)^{1/4}\mathrm{Ai}(f_n(z)) - \left(\frac{z+1}{z-1}f_n(z)\right)^{-1/4}\mathrm{Ai}'(f_n(z))\right\},$$

the error term being uniform in $C_{2,\delta_0}$ for some $\delta_0 > 0$. This captures the case $u \geq 0$ in the statement of Lemma 7.2. The case $u < 0$ would be captured by a similar estimate valid in the region $C_{1,\delta}$. To explain (7.5), the equilibrium measure for the Gaussian distribution is given by the well-known Wigner semicircle

$$\Psi(y) = \frac{2}{\pi}(1-y^2)^{1/2}$$

and

$$f_n(z) = -(-n)^{2/3}\left(\frac{3\pi}{2}\right)^{2/3}\left(\int_1^z \Psi(y)\,dy\right)^{2/3}$$
$$= -\left(-\frac{3n}{2}\right)^{2/3}(z\sqrt{1-z^2} - \arccos z)^{2/3}.$$

Setting $z = 1 + x$ for small $x \geq 0$, one computes

$$f_n(1+x) = 2n^{2/3}x\left(1 + \frac{x}{10} + \cdots\right)$$

and

$$\left(\frac{z+1}{z-1}f_n(z)\right)^{1/4}\bigg|_{z=1+x} = \sqrt{2}\,n^{1/6}\left(1 + \frac{3}{20}x + \cdots\right).$$

Defining $x$ such that

$$\sqrt{2n+1} + \frac{u}{\sqrt{2}n^{1/6}} = \sqrt{2(n-k)}(1+x),$$

one computes, for $k = Mn^{1/3} \log n$ and $|M| \leq M_0$,

$$
\begin{aligned}
(7.6) \quad x &= \frac{u}{2n^{2/3}} + \frac{2k+1}{4n} + \frac{ku}{4n^{5/3}} + \cdots \\
&= \frac{M \log n + u}{2n^{2/3}} + \frac{1}{4n} + \frac{Mu \log n}{4n^{4/3}} + \cdots.
\end{aligned}
$$

Thus for $x$ behaving as (7.6) and for $k = Mn^{1/3} \log n$, we deduce from the formulae above

$$
\begin{aligned}
f_{n-k}(z)|_{z=1+x} &= 2(n-k)^{2/3} x \left(1 + \frac{x}{10} + \cdots\right) \\
&= 2n^{2/3} x \left(1 - \frac{2k}{3n} + \cdots\right)\left(1 + \frac{x}{10} + \cdots\right) \\
&= M \log n + u + \frac{1}{2n^{1/3}} + O\left(\frac{M \log n}{n^{2/3}}\right),
\end{aligned}
$$

$$
\begin{aligned}
\left(\frac{z+1}{z-1} f_{n-k}(z)\right)^{1/4}\bigg|_{z=1+x} &= \sqrt{2}(n-k)^{1/6}\left(1 + \frac{3}{20} x + \cdots\right) \\
&= \sqrt{2} n^{1/6}\left(1 - \frac{k}{6n} + \cdots\right)\left(1 + \frac{3}{20} x + \cdots\right) \\
&= \sqrt{2}\, n^{1/6}\left(1 + O\left(\frac{M \log n}{n^{2/3}}\right)\right)
\end{aligned}
$$

and

$$
\frac{1}{(n-k)^{1/4}} = \frac{1}{n^{1/4}}\left(1 + \frac{k}{4n} + \cdots\right) = \frac{1}{n^{1/4}}\left(1 + O\left(\frac{M \log n}{n^{2/3}}\right)\right).
$$

Using the asymptotics (7.5), one computes for $k = Mn^{1/3} \log n$ and $|M| \leq M_0$,

$$
\begin{aligned}
&\varphi_{n-k}\left(\sqrt{2n+1} + \frac{u}{\sqrt{2} n^{1/6}}\right) \\
&= \varphi_{n-k}(\sqrt{2(n-k)}(1+x)) \\
&= \frac{1 + O(1/(n-k))}{(2(n-k))^{1/4}}\Bigg\{ \left(\frac{z+1}{z-1} f_{n-k}(z)\right)^{1/4} \mathrm{Ai}(f_{n-k}(z)) \\
&\qquad\qquad - \left(\frac{z+1}{z-1} f_{n-k}(z)\right)^{-1/4} \mathrm{Ai}'(f_{n-k}(z))\Bigg\}\bigg|_{z=1+x} \\
&= \frac{1 + O(1/n)}{(2n)^{1/4}}\left(1 + O\left(\frac{M \log n}{n^{2/3}}\right)\right) \\
&\quad \times \left\{\sqrt{2} n^{1/6}\left(1 + O\left(\frac{M \log n}{n^{2/3}}\right)\right)\right.
\end{aligned}
$$

$$\times \operatorname{Ai}\left(M\log n + u + \frac{1}{2n^{1/3}} + O\left(\frac{M\log n}{n^{2/3}}\right)\right)$$

$$-\frac{1}{\sqrt{2}n^{1/6}}\left(1 + O\left(\frac{M\log n}{n^{2/3}}\right)\right)$$

$$\times \operatorname{Ai}'\left(M\log n + u + \frac{1}{2n^{1/3}} + O\left(\frac{\log n}{n^{2/3}}\right)\right)\Bigg\}$$

$$= \frac{2^{1/4}(1 + O(M\log n/n^{2/3}))}{n^{1/12}}$$

$$\times \Bigg\{\left(1 + O\left(\frac{M\log n}{n^{2/3}}\right)\right)$$

$$\times \left[\operatorname{Ai}(M\log n + u) + \frac{\operatorname{Ai}'(M\log n + u)}{2n^{1/3}} + O\left(\frac{M\log n}{n^{2/3}}\right)\right]$$

$$-\left(1 + O\left(\frac{M\log n}{n^{2/3}}\right)\right)\left[\frac{\operatorname{Ai}'(M\log n + u)}{2n^{1/3}} + O\left(\frac{M\log n}{n^{2/3}}\right)\right]\Bigg\}$$

$$= \frac{2^{1/4}}{n^{1/12}}\operatorname{Ai}(u + M\log n)\left[1 + O\left(\frac{M\log n}{n^{2/3}}\right)\right]$$

$$= \frac{2^{1/4}}{n^{1/12}}\operatorname{Ai}\left(u + \frac{k}{n^{1/3}}\right)\left[1 + O\left(\frac{M\log n}{n^{2/3}}\right)\right],$$

ending the proof of Lemma 7.2. □

PROOF OF PROPOSITION 7.1. As a first step, in a recent paper Krasikov [10] shows the following inequality for $k \geq 6$ and for a universal constant $c$:

$$(7.7) \qquad \max_x |H_k(x)|e^{-x^2/2} \leq c\frac{k!}{(k/2)!}\left(\frac{2k\sqrt{4k-2}}{k^{1/6}\sqrt{8k^2 - 8k + 3}}\right)^{1/2}.$$

From Stirling's formula $n! = \sqrt{2\pi n}\, n^n e^{-n}(1 + O(1/n))$, it follows that

$$\frac{\sqrt{k!}}{2^{k/2}(k/2)!} = \left(\frac{2}{\pi k}\right)^{1/4}\left(1 + O\left(\frac{1}{k}\right)\right).$$

This estimate, estimate (7.7) and formula (7.2) show that, for $k \geq$ some fixed $k_0$ and some constant $c'$,

$$(7.8) \qquad \max_x |\varphi_k(x)| = \frac{1}{2^{k/2}\sqrt{k!}\pi^{1/4}}\max_x |H_k(x)|e^{-x^2/2} \leq \frac{c'}{k^{1/12}}.$$

Using both estimates (7.4) and (7.8) in (7.1), one computes for $t \geq s$, taking into account substitution $\mathcal{S}_1$, as in (7.3),

$$\left|K_{t,s}^{H,n}(x,y)\,dy_{|\mathcal{S}_1}\right.$$

$$
\begin{aligned}
&- \sum_{k=0}^{[Mn^{1/3}\log n]} e^{-(k/n^{1/3})(t-s)} \mathrm{Ai}\left(u + \frac{k}{n^{1/3}}\right) \\
&\qquad\qquad \times \mathrm{Ai}\left(v + \frac{k}{n^{1/3}}\right) \left(\frac{2^{1/4}}{n^{1/12}}\right)^2 \frac{dv}{\sqrt{2}n^{1/6}} \Bigg| \\
&\leq \Bigg| \sum_{k=0}^{n} e^{-(k/n^{1/3})(t-s)} \varphi_{n-k}\left(\sqrt{2n+1} + \frac{u}{\sqrt{2}n^{1/6}}\right) \\
&\qquad\qquad \times \varphi_{n-k}\left(\sqrt{2n+1} + \frac{v}{\sqrt{2}n^{1/6}}\right) \\
&\quad - \sum_{k=0}^{[Mn^{1/3}\log n]} e^{-(k/n^{1/3})(t-s)} \frac{2^{1/4}}{n^{1/12}} \mathrm{Ai}\left(u + \frac{k}{n^{1/3}}\right)(1 + E_u(k,n)) \\
&\qquad\qquad \times \frac{2^{1/4}}{n^{1/12}} \mathrm{Ai}\left(v + \frac{k}{n^{1/3}}\right)(1 + E_v(k,n)) \Bigg| \frac{dv}{\sqrt{2}n^{1/6}} \\
&\quad + \sum_{k=0}^{[Mn^{1/3}\log n]} e^{-(k/n^{1/3})(t-s)} \frac{2^{1/4}}{n^{1/6}} \left| \mathrm{Ai}\left(u + \frac{k}{n^{1/3}}\right) \mathrm{Ai}\left(v + \frac{k}{n^{1/3}}\right) \right| \frac{dv}{\sqrt{2}n^{1/6}} \\
&\qquad\qquad \times |1 - (1 + E_u(k,n))(1 + E_v(k,n))| \\
&\leq \Bigg| \sum_{k=[Mn^{1/3}\log n]+1}^{n} e^{-(k/n^{1/3})(t-s)} \\
&\qquad\qquad \times \varphi_{n-k}\left(\sqrt{2n+1} + \frac{u}{\sqrt{2}n^{1/6}}\right) \\
&\qquad\qquad \times \varphi_{n-k}\left(\sqrt{2n+1} + \frac{v}{\sqrt{2}n^{1/6}}\right) \Bigg| \frac{dv}{\sqrt{2}n^{1/6}} \\
&\quad + \frac{c''}{n^{1/6}}(Mn^{1/3}\log n)\frac{\log n}{n^{2/3}}(1+o(1))\frac{dv}{\sqrt{2}n^{1/6}} \\
&\leq \frac{dv}{\sqrt{2}n^{1/6}}\Bigg( c'^2 \sum_{m=0}^{n-1-[Mn^{1/3}\log n]} e^{-(n-m)(t-s)/n^{1/3}} \frac{1}{m^{1/6}} \\
&\qquad\qquad + c''M \frac{(\log n)^2}{n^{1/2}}(1+o(1)) \Bigg),
\end{aligned}
$$
(7.9)

where $c''$ is determined by the maximum of the Airy function $\mathrm{Ai}(z)$ on the semi-infinite interval $(0,\infty)$.

Setting $n' = n - 1 - [Mn^{1/3}\log n]$, the sum in the last expression is estimated as follows:

$$\frac{1}{\sqrt{2}n^{1/6}}\sum_{m=1}^{n'}\frac{1}{m^{1/6}}e^{-(n-m)(t-s)/n^{1/3}}$$

$$\leq \frac{1}{\sqrt{2}n^{1/6}}\left(\sum_{0}^{\ell-1}e^{-(n-m)(t-s)/n^{1/3}} + \frac{1}{\ell^{1/6}}\sum_{\ell}^{n'}e^{-(n-m)(t-s)/n^{1/3}}\right)$$

$$= \frac{1}{\sqrt{2}n^{1/6}}\frac{1}{1-e^{-(t-s)/n^{1/3}}}[e^{-(n-\ell+1)(t-s)/n^{1/3}} - e^{-(n+1)(t-s)/n^{1/3}}$$

$$+ \ell^{-1/6}(e^{-(n-n')(t-s)/n^{1/3}} - e^{-(n-\ell+1)(t-s)/n^{1/3}})]$$

$$\simeq \frac{1}{\sqrt{2}n^{1/6}}\frac{n^{1/3}}{t-s}[e^{-(n-\ell+1)(t-s)/n^{1/3}} - e^{-(n+1)(t-s)/n^{1/3}}$$

$$+ \ell^{-1/6}(e^{-(n-n')(t-s)/n^{1/3}}(*) - e^{-(n-\ell+1)(t-s)/n^{1/3}})].$$

Picking $\ell = O(n/2)$, all terms above tend to $0$ exponentially fast, except the term $(*)$, which requires some attention. Choosing $n' = n - [Mn^{1/3}\log n] - 1$, so that $n - n' = O(n^{1/3}\log n)$, the coefficient of that term is $O(n^{-M(t-s)})$. Therefore that term gets small, when $n \to \infty$ and $\ell = O(n/2)$,

$$\frac{1}{\sqrt{2}n^{1/6}}\frac{n^{1/3}}{\ell^{1/6}}O(n^{-M(t-s)}) \to 0.$$

The proof is ended by observing that the second term in the first difference of (7.9) is a Riemann sum converging to the extended Airy kernel, that is,

$$\sum_{k=0}^{[Mn^{1/3}\log n]} e^{-(k/n^{1/3})(t-s)}\text{Ai}\left(u+\frac{k}{n^{1/3}}\right)\text{Ai}\left(v+\frac{k}{n^{1/3}}\right)\left(\frac{2^{1/4}}{n^{1/12}}\right)^2\frac{dv}{\sqrt{2}n^{1/6}}$$

$$= dv\int_0^\infty e^{-z(t-s)}\text{Ai}(u+z)\text{Ai}(v+z)\,dz + O\left(\frac{M\log n}{n^{1/3}}\right).$$

This establishes the convergence for $t \geq s$.

For $t < s$, one computes, again using in (7.1) the estimates (7.4) of Lemma 7.2 and (7.8),

$$\left|K_{t,s}^{H,n}(x,y)\,dy_{|S_1}\right.$$

$$\left. + \sum_{k=0}^{[M_0 n^{1/3}\log n]} e^{-(k/n^{1/3})(s-t)}\text{Ai}\left(u-\frac{k}{n^{1/3}}\right)\text{Ai}\left(v-\frac{k}{n^{1/3}}\right)\left(\frac{2^{1/4}}{n^{1/12}}\right)^2\frac{dv}{\sqrt{2}n^{1/6}}\right|$$

$$\leq \left| \sum_{k=0}^{\infty} e^{-(k/n^{1/3})(s-t)} \varphi_{n+k}\left(\sqrt{2n+1} + \frac{u}{\sqrt{2}n^{1/6}}\right) \right.$$

$$\times \varphi_{n+k}\left(\sqrt{2n+1} + \frac{v}{\sqrt{2}n^{1/6}}\right)$$

$$- \sum_{k=0}^{[M_0 n^{1/3} \log n]} e^{-(k/n^{1/3})(s-t)} \frac{2^{1/4}}{n^{1/12}} \mathrm{Ai}\left(u - \frac{k}{n^{1/3}}\right)(1 + E(k,n))$$

$$\times \left. \frac{2^{1/4}}{n^{1/12}} \mathrm{Ai}\left(v - \frac{k}{n^{1/3}}\right)(1 + E(k,n)) \right| \frac{dv}{\sqrt{2}n^{1/6}}$$

$$+ \sum_{k=0}^{[Mn^{1/3} \log n]} e^{-(k/n^{1/3})(t-s)} \frac{2^{1/4}}{n^{1/6}} \left| \mathrm{Ai}\left(u + \frac{k}{n^{1/3}}\right) \mathrm{Ai}\left(v + \frac{k}{n^{1/3}}\right) \right| \frac{dv}{\sqrt{2}n^{1/6}}$$

(7.10) $$\times |1 - (1 + E_u(k,n))(1 + E_v(k,n))|$$

$$\leq \left| c'' M \frac{(\log n)^2}{n^{1/2}}(1 + o(1)) \right.$$

$$+ \sum_{k=[Mn^{1/3}\log n]+1}^{\infty} e^{-(k/n^{1/3})(s-t)}$$

$$\times \varphi_{n+k}\left(\sqrt{2n+1} + \frac{u}{\sqrt{2}n^{1/6}}\right)$$

$$\left. \times \varphi_{n+k}\left(\sqrt{2n+1} + \frac{v}{\sqrt{2}n^{1/6}}\right) \right| \frac{dv}{\sqrt{2}n^{1/6}}$$

$$\leq \frac{dv}{\sqrt{2}n^{1/6}} \left( c'^2 \sum_{k=[Mn^{1/3}\log n]+1}^{\infty} \frac{e^{-(k/n^{1/3})(s-t)}}{(n+k)^{1/6}} + c''M \frac{(\log n)^2}{n^{1/2}}(1+o(1)) \right)$$

$$\leq \frac{dv}{\sqrt{2}n^{1/6}} \left( c'^2 n^{-M(s-t)-1/6} \sum_{m=1}^{\infty} e^{-(m/n^{1/3})(s-t)} \right.$$

$$\left. + c''M \frac{(\log n)^2}{n^{1/2}}(1+o(1)) \right).$$

The rest of the proof goes the same way as before, ending the proof of Proposition 7.1. $\square$

7.3. *The extended Hermite kernel tends to the extended Sine kernel.* Given the substitution

(7.11) $$\mathcal{S}_2 := \left\{ t \mapsto \frac{\pi^2 t}{2n}, s \mapsto \frac{\pi^2 s}{2n}, x \mapsto \frac{\pi u}{\sqrt{2n}}, y \mapsto \frac{\pi v}{\sqrt{2n}} \right\},$$

we have:

PROPOSITION 7.3. *The extended Hermite kernel tends to the extended Sine kernel, when $n \to \infty$, uniformly for $u, v \in$ compact subsets $\subset \mathbb{R}$:*

$$K_{t,s}^{H,n}(x,y)\, dy_{|\mathcal{S}_2} \longrightarrow e^{-(\pi^2/2)(t-s)} K_{t,s}^S(u,v)\, dv.$$

SKETCH OF PROOF. From [13], page 198, it follows that for $|x| \leq M$,

$$\varphi_k(x) = \frac{1}{2^{k/2}\sqrt{k!}\pi^{1/4}} e^{-x^2/2} H_k(x)$$

$$= \frac{\sqrt{k!}}{2^{k/2}(k/2)!\pi^{1/4}} \left( \cos\left( x\sqrt{2k+1} - \frac{k\pi}{2} \right) \right.$$

$$\left. + \frac{x^3}{6\sqrt{2k+1}} \sin\left( x\sqrt{2k+1} - \frac{k\pi}{2} \right) + O(k^{-1}) \right)$$

$$= \left(\frac{2}{k\pi^2}\right)^{1/4} \left( \cos\left(\sqrt{2k+1}\, x - \frac{k\pi}{2}\right) + O\left(\frac{1}{\sqrt{k}}\right) \right).$$

Using the substitution $\mathcal{S}_2$ as in (7.11), one computes in (7.1) for $t > s$,

(7.12) 
$$\begin{aligned}
&\left| K_{t,s}^{H,n}(x,y)\, dy_{|\mathcal{S}} \right. \\
&\quad \left. - e^{-\pi^2(t-s)/2} \frac{\pi}{\sqrt{2n}} \sum_{k=\ell+1}^{n-1} e^{\pi^2(t-s)k/2n} \varphi_k\left(\frac{\pi u}{\sqrt{2n}}\right) \varphi_k\left(\frac{\pi v}{\sqrt{2n}}\right) dv \right| \\
&= \left| e^{-\pi^2(t-s)/2} \frac{\pi}{\sqrt{2n}} \sum_{k=0}^{\ell} e^{\pi^2(t-s)k/2n} \varphi_k\left(\frac{\pi u}{\sqrt{2n}}\right) \varphi_k\left(\frac{\pi v}{\sqrt{2n}}\right) dv \right| \\
&\leq e^{-\pi^2(t-s)/2} \frac{\pi}{\sqrt{2n}} \left( \sum_{k=0}^{\alpha} + \sum_{k=\alpha+1}^{\ell} \right) \\
&\quad \times e^{\pi^2(t-s)k/2n} \left| \varphi_k\left(\frac{\pi u}{\sqrt{2n}}\right) \varphi_k\left(\frac{\pi v}{\sqrt{2n}}\right) dv \right| \\
&= (\mathrm{I}) + (\mathrm{II}),
\end{aligned}$$

where $\alpha > 0$ is the minimal integer above which Krasikov's estimate (7.8) holds. But then expression (I) tends to 0 exponentially, when $n$ tends to $\infty$, and (II) is estimated as follows:

$$(\mathrm{II}) \leq \frac{c'^2 \pi}{\sqrt{2n}} \sum_{k=\alpha+1}^{\ell} \frac{1}{k^{1/6}} \leq \frac{c'' \pi}{\sqrt{2n}} \ell^{5/6},$$

which tends to 0 for $\ell, n \to \infty$ such that $\ell^{5/6}/n^{1/2} \to 0$. Then the sum appearing at the first line of (7.12) can be estimated as follows:

$$e^{-\pi^2(t-s)/2}\frac{\pi}{\sqrt{2n}}\sum_{k=\ell+1}^{n-1} e^{\pi^2(t-s)k/(2n)}\varphi_k\left(\frac{\pi u}{\sqrt{2n}}\right)\varphi_k\left(\frac{\pi v}{\sqrt{2n}}\right)dv$$

$$= e^{-\pi^2(t-s)/2}\sum_{k=\ell+1}^{n-1}\frac{1}{n}e^{\pi^2 k(t-s)k/(2n)}\left(\frac{n}{k}\right)^{1/2}$$

$$\times\left[\cos\left(\pi u\sqrt{\frac{k}{n}+\frac{1}{2n}}-\frac{n\pi}{2}\frac{k}{n}\right)\right.$$

$$\left.\times\cos\left(\pi v\sqrt{\frac{k}{n}+\frac{1}{2n}}-\frac{n\pi}{2}\frac{k}{n}\right)+O\left(\frac{1}{\sqrt{\ell}}\right)\right]dv$$

$$= e^{-\pi^2(t-s)/2}\sum_{k=\ell+1}^{n-1}\frac{1}{2n}e^{\pi^2 k(t-s)/(2n)}\left(\frac{n}{k}\right)^{1/2}$$

$$\times\left[\cos\left(\pi(u-v)\sqrt{\frac{k}{n}+\frac{1}{2n}}\right)\right.$$

$$\left.+\cos\left(\pi(u+v)\sqrt{\frac{k}{n}+\frac{1}{2n}}-n\pi\frac{k}{n}\right)\right.$$

$$\left.+O\left(\frac{1}{\sqrt{\ell}}\right)\right]dv$$

$$\simeq\int_{\ell/n}^{1}\frac{dx}{2\sqrt{x}}e^{\pi^2(x-1)(t-s)/2}$$

$$\times\left[\cos(\pi(u-v)\sqrt{x})\right.$$

$$\left.+\cos(\pi(u+v)\sqrt{x}-n\pi x)+O\left(\frac{1}{\sqrt{\ell}}\right)+O\left(\frac{1}{\sqrt{n}}\right)\right]dv$$

$$=\frac{1}{\pi}\int_{\pi\sqrt{\ell/n}}^{\pi}dz\,e^{-(\pi^2-z^2)(t-s)/2}$$

$$\times\left[\cos(z(u-v))+\cos(z(u+v)-nz^2/\pi)\right.$$

$$\left.+O\left(\frac{1}{\sqrt{\ell}}\right)+O\left(\frac{1}{\sqrt{n}}\right)\right]dv(*)$$

$$\to \frac{1}{\pi} e^{-\pi^2(t-s)/2} \int_0^\pi dz\, e^{(z^2/2)(t-s)} \cos(z(u-v))\, dv$$
$$= e^{-\pi^2(t-s)/2} K_{t,s}^S(u,v)$$

when $\ell$ and $n$ tend to $\infty$. The integral involving the second cosine in the expression $(*)$ above is an oscillatory integral and thus tends to zero faster than any power of $n$ by the Riemann–Lebesgue lemma. The case $t \leq s$ proceeds along similar lines, establishing Proposition 7.3. $\square$

7.4. *Convergence of Fredholm determinants and their derivatives.* Here we give a schematic argument, based on the customary formula $\log\det(I - K) = \operatorname{Tr}\log(I - K)$, where $K$ is a kernel restricted to a disjoint union of intervals $E$ and where $K + \delta K$ tends to $K$. Then, from

$$(7.13) \quad \begin{aligned} \det(I - K - \delta K) &= \det(I - K)\det(I - (I - K)^{-1}\delta K) \\ &= \det(I - K)(1 - \operatorname{Tr}(I - K)^{-1}\delta K + o(\delta K)), \end{aligned}$$

one sees that $\det(I - K - \delta K)$ tends to $\det(I - K)$, when $\delta K$ tends to $0$. Also, given $p_1, \ldots, p_{2r}$ the endpoints of the set $E$ (see [7, 15]),

$$\frac{\partial}{\partial p_k} \log\det(I - K) = (-1)^{k-1} K(I - K)^{-1}(p_k, p_k),$$

where here "=" means "kernel of," evaluated at $(p_k, p_k)$ and so

$$\frac{\partial}{\partial p_k} \log\det(I - (K + \delta K))$$
$$= (-1)^{k-1}(K + \delta K)(I - K - \delta K)^{-1}(p_k, p_k)$$
$$= (-1)^{k-1}[K(I - K)^{-1}(p_k, p_k)$$
$$\qquad + (\delta K(I - K)^{-1}$$
$$\qquad\qquad + K(I - K)^{-1}\delta K(I - K)^{-1})(p_k, p_k) + o(\delta K)]$$
$$= \frac{\partial}{\partial p_k} \log\det(I - K) + O(\delta K).$$

Since by Propositions 7.1 and 7.3 the extended Hermite kernel converges to the extended Airy and Sine kernels, this argument shows the convergence of the corresponding Fredholm determinants and their first derivatives with respect to the end points of $E$. In a similar fashion one proves the result for higher derivatives.

7.5. *An a priori asymptotic expansion for the joint Airy probability.* The proof of Theorem 1.6 in Section 6 was based on an a priori asymptotic expansion for the ratio below in $1/t$ for large $t = t_2 - t_1$. This can be found in [18] and proceeds as follows:

$$
\begin{aligned}
\frac{P(A(t_1) \leq u, A(t_2) \leq v)}{P(A(t_1) \leq u) P(A(t_2) \leq v)} &= \frac{\det(I - (\hat{K}^A_{t_i t_j})_{1 \leq i,j \leq 2})}{\det(I - \hat{K}^A_{t_1 t_1}) \det(I - \hat{K}^A_{t_2 t_2})} \\
&= \det\left(I - \begin{pmatrix} 0 & \mathcal{K}_{12} \\ \mathcal{K}_{21} & 0 \end{pmatrix}\right) \\
&= \det(I - \mathcal{K}_{12} \mathcal{K}_{21}) \\
&= 1 + \sum_{i \geq 1} \frac{f_i(u,v)}{t^i},
\end{aligned}
\tag{7.14}
$$

where

$$\mathcal{K}_{12} := (I - \chi_{[u,\infty)}(x) K^A_{00}(x,y) \chi_{[u,\infty)}(y))^{-1} \chi_{[u,\infty)}(x) K^A_{0,t}(x,y) \chi_{[v,\infty)}(y),$$

$$\mathcal{K}_{21} := (I - \chi_{[v,\infty)}(x) K^A_{00}(x,y) \chi_{[v,\infty)}(y))^{-1} \chi_{[v,\infty)}(x) K^A_{t,0}(x,y) \chi_{[u,\infty)}(y),$$

with

$$K^A_{00}(x,y) = \int_0^\infty \operatorname{Ai}(x+z) \operatorname{Ai}(y+z) \, dz,$$

$$K^A_{t,0}(x,y) = \int_0^\infty e^{-zt} \operatorname{Ai}(x+z) \operatorname{Ai}(y+z) \, dz = O(1/t),$$

$$K^A_{0,t}(x,y) = -\int_{-\infty}^0 e^{zt} \operatorname{Ai}(x+z) \operatorname{Ai}(y+z) \, dz = O(1/t).$$

## APPENDIX: REMARK ABOUT THE "NONEXPLOSION" CONJECTURE

To discuss the conjecture (1.21), consider the Dyson Brownian motion $(\lambda_1(t), \ldots, \lambda_n(t))$ and the corresponding Ornstein–Uhlenbeck process on the matrix $B$. Then, using the change of variables

$$M_i = \frac{B_i}{\sqrt{(1-c^2)/2}},$$

and further $M_2 \mapsto M := M_2 - cM_1$ in the $M_2$-integrals below and noting that $\max(\operatorname{spec} M_1) \leq -z$ and $\max(\operatorname{spec} M_2) \geq a$ imply $\max(\operatorname{spec}(M_2 - cM_1)) \geq a + cz$, we have for the conditional probability, the following inequality:

$$P(\lambda_n(t) \geq a \mid \lambda_n(0) \leq -z)$$

$$= \frac{\int_{\max(\operatorname{spec} M_1) \leq -z} dM_1 \, e^{-(1-c^2) \operatorname{Tr} M_1^2 / 2} \int_{\max(\operatorname{spec} M_2) \geq a} dM_2 \, e^{-\operatorname{Tr}(M_2 - cM_1)^2 / 2}}{\int_{\max(\operatorname{spec} M_1) \leq -z} dM_1 \, e^{-(1-c^2) \operatorname{Tr} M_1^2 / 2} \int_{M_2 \in \mathcal{H}_n} dM_2 \, e^{-\operatorname{Tr}(M_2 - cM_1)^2 / 2}}$$

$$\leq \frac{\int_{\max(\operatorname{spec} M_1) \leq -z} dM_1 \, e^{-(1-c^2) \operatorname{Tr} M_1^2 / 2} \int_{\max(\operatorname{spec} M) \geq a + cz} dM \, e^{-\operatorname{Tr} M^2 / 2}}{\int_{\max(\operatorname{spec} M_1) \leq -z} dM_1 \, e^{-(1-c^2) \operatorname{Tr} M_1^2 / 2} \int_{M \in \mathcal{H}_n} dM \, e^{-\operatorname{Tr} M^2 / 2}}$$

$$= P(\lambda_n(t) \geq a + cz),$$

implying

$$\lim_{z \to \infty} P(\lambda_n(t) \geq a \mid \lambda_n(0) \leq -z) = 0,$$

and a fortiori,

$$\lim_{z \to \infty} P(\lambda_n(t) \geq x + z \mid \lambda_n(0) \leq -z) = 0.$$

It is unclear why the limit (6.5) remains valid when $n \to \infty$, using the Airy scaling (1.4).

Department of Mathematics  
Brandeis University  
Waltham, Massachusetts 02454  
USA  
e-mail: adler@brandeis.edu

Department of Mathematics  
Université de Louvain  
1348 Louvain-la-Neuve  
Belgium  
and  
Brandeis University  
Waltham, Massachusetts 02454  
USA  
e-mail: vanmoerbeke@math.ucl.ac.be  
e-mail: vanmoerbeke@brandeis.edu